\documentclass[11pt]{article}

\input xy

\xyoption{all}

\usepackage{amssymb,amsbsy,amsmath,graphicx,epsfig,times}

\newtheorem{thm}{Theorem}[section]
\newtheorem{lem}[thm]{Lemma}
\newtheorem{prop}[thm]{Proposition}
\newtheorem{cor}[thm]{Corollary}

\newtheorem{dfn}[thm]{Definition}

\renewcommand{\bf}[1]{\mathbf{#1}}
\renewcommand{\rm}[1]{\mathrm{#1}}
\renewcommand{\cal}[1]{\mathcal{#1}}

\newcommand{\bbC}{\mathbb{C}}

\newcommand{\bbN}{\mathbb{N}}

\newcommand{\bbR}{\mathbb{R}}
\newcommand{\bbT}{\mathbb{T}}
\newcommand{\bbZ}{\mathbb{Z}}

\newcommand{\bfW}{\mathbf{W}}
\newcommand{\bfX}{\mathbf{X}}
\newcommand{\bfY}{\mathbf{Y}}
\newcommand{\bfZ}{\mathbf{Z}}

\newcommand{\sfC}{\mathsf{C}}

\newcommand{\sfE}{\mathsf{E}}

\renewcommand{\d}{\mathrm{d}}

\newcommand{\A}{\mathcal{A}}

\newcommand{\F}{\mathcal{F}}

\newcommand{\frH}{\mathfrak{H}}

\newcommand{\frg}{\mathfrak{g}}
\newcommand{\frh}{\mathfrak{h}}

\newcommand{\G}{\Gamma}
\renewcommand{\L}{\Lambda}

\renewcommand{\S}{\Sigma}

\renewcommand{\a}{\alpha}

\renewcommand{\l}{\lambda}
\newcommand{\s}{\sigma}
\renewcommand{\phi}{\varphi}

\newcommand{\id}{\mathrm{id}}

\renewcommand{\hat}[1]{\widehat{#1}}

\renewcommand{\to}{\longrightarrow}
\newcommand{\fin}{\nolinebreak\hspace{\stretch{1}}$\lhd$}
\newcommand{\qed}{\nolinebreak\hspace{\stretch{1}}$\Box$}

\renewcommand{\t}[1]{\tilde{#1}}
\newcommand{\actson}{\curvearrowright}\newcommand{\barint}{-\!\!\!\!\!\!\int}
\newcommand{\img}{\rm{img}\,}

\newcommand{\PET}{\rm{PET}}
\newcommand{\Lat}{\rm{Lat}\,}
\newcommand{\onto}{\twoheadrightarrow}
\newcommand{\ldeg}{\rm{ldeg}\,}
\newcommand{\LT}{\rm{LT}}
\newcommand{\wt}{\rm{wt}\,}
\newcommand{\Wt}{\rm{Wt}\,}
\renewcommand{\vec}[1]{\stackrel{\rightarrow}{#1}}

\begin{document}

\title{Equidistribution of joinings under off-diagonal polynomial flows of nilpotent Lie groups}
\author{Tim Austin}

\date{}

\maketitle

\begin{abstract}
Let $G$ be a connected nilpotent Lie group.  Given probability-preserving $G$-actions $(X_i,\S_i,\mu_i,u_i)$, $i=0,1,\ldots,k$, and also polynomial maps $\phi_i:\bbR\to G$, $i=1,\ldots,k$, we consider the trajectory of a joining $\l$ of the systems $(X_i,\S_i,\mu_i,u_i)$ under the `off-diagonal' flow
\[(t,(x_0,x_1,x_2,\ldots,x_k))\mapsto (x_0,u_1^{\phi_1(t)}x_1,u_2^{\phi_2(t)}x_2,\ldots,u_k^{\phi_k(t)}x_k).\]
It is proved that any joining $\l$ is equidistributed under this flow with respect to some limit joining $\l'$.  This is deduced from the stronger fact of norm convergence for a system of multiple ergodic averages, related to those arising in Furstenberg's approach to the study of multiple recurrence.  It is also shown that the limit joining $\l'$ is invariant under the subgroup of $G^{k+1}$ generated by the image of the off-diagonal flow, in addition to the diagonal subgroup.
\end{abstract}

\parskip 0pt

\tableofcontents

\parskip 7pt

\parindent 0pt

\section{Introduction}

This paper is set among jointly measurable probability-preserving actions of a connected nilpotent Lie group $G$.  We will assume in addition that $G$ is simply connected; it will be clear from the statements of our main results that by ascending to the universal cover this incurs no real loss of generality.

Suppose that $u_i:G\actson (X_i,\S_i,\mu_i)$ for $i=0,1,\ldots,k$ is a tuple of such actions and that $\l$ is a joining of them.  This means that $\l$ is a coupling of the measures $\mu_i$ on the product space $\prod_iX_i$, and that it is invariant under the \textbf{diagonal transformation}
\[u_\Delta^g := u_0^g\times u_1^g\times \cdots\times u_k^g\]
for every $g \in G$.

Taking the $G$-actions on each coordinate separately, the $u_i$ together define a jointly measurable action $u_\times$ of the whole Cartesian power $G^{k+1}$ on $\prod_i X_i$ according to
\[u_\times^{(g_0,g_1,\ldots,g_k)}:= u_0^{g_0}\times u_1^{g_1}\times \cdots\times u_k^{g_k}.\]
In these terms $u_\Delta$ may be identified with the restriction of $u_\times$ to the diagonal subgroup
\[G^{\Delta (k+1)}:= \{(g,g,\ldots,g):\ g\in G\}\leq G^{k+1}.\]

An arbitrary joining $\l$ need not be $u_\times$-invariant.  However, the main result of this paper implies that for any one-parameter subgroup $\bbR\to G^{k+1}$, the trajectory of $\l$ under the $u_\times$-action of that subgroup must equidistribute with respect to some new joining $\l'$ that is also invariant under that subgroup.  Moreover, this statement generalizes to averages over the trajectory of any map $\bbR\to G^{k+1}$ that is `polynomial' in the sense that repeated group-valued differencing leads to the trivial map (precise definitions are recalled in Section~\ref{sec:poly}).  The full result is the following.

\begin{thm}\label{thm:main}
If $(X_i,\S_i,\mu_i,u_i)$, $0 \leq i \leq k$, and $\l$ are as above, and $\phi_i:\bbR\to G$ for $1 \leq i \leq k$ are polynomial maps satisfying $\phi_i(0) = e$ (the identity of $G$), then the averaged measures
\[\l_T := \barint_0^T (\id_{X_0} \times u_1^{\phi_1(t)}\times \cdots\times u_k^{\phi_k(t)})_\ast\l\,\d t\]
converge in the coupling topology as $T\to\infty$ to some joining $\l'$ of the systems $(X_i,\S_i,\mu_i,u_i)$ which is invariant under the restriction of $u_\times$ to the subgroup
\[\langle G^{\Delta (k+1)}\cup \{(e,\phi_1(t),\ldots,\phi_k(t)):\ t\in\bbR\}\rangle.\]
\end{thm}

Here we have used the standard analyst's notation $\barint_a^b := \frac{1}{b-a}\int_a^b$, and we write $\langle S\rangle$ for the smallest closed subgroup of $G$ containing $S$.

\textbf{Remark}\quad If $t$ is such that $\phi_i(t) \neq e$ then the individual measures
\[(\id_{X_0} \times u_1^{\phi_1(t)}\times \cdots\times u_k^{\phi_k(t)})_\ast\l\]
may \emph{not} be joinings of the original actions.  As measures they are still couplings of the $\mu_i$, but the invariance of $\l$ under the diagonal subgroup has been replaced with invariance under its conjugate
\[(e,\phi_1(t),\ldots,\phi_k(t))\cdot G^{\Delta (k+1)}\cdot (e,\phi_1(t),\ldots,\phi_k(t))^{-1}.\]
Thus a non-trivial part of the conclusion of Theorem~\ref{thm:main} is that the smoothing effect of averaging over $t$ recovers the invariance under $G^{\Delta (k+1)}$ (and  likewise under all of these conjugates). \fin

Convergence $\l_T\to \l'$ in the coupling topology, as in Theorem~\ref{thm:main}, asserts that
\[\int_{X_0\times X_1\times \cdots\times X_k}f_0\otimes f_1\otimes \cdots \otimes f_k\,\d\l_T\to \int_{X_0\times X_1\times \cdots\times X_k}f_0\otimes f_1\otimes \cdots \otimes f_k\,\d\l'\]
for any choice of $f_0\in L^\infty(\mu_0)$, $f_1 \in L^\infty(\mu_1)$, \ldots, $f_k \in L^\infty(\mu_k)$.  Informally, it is a variant of weak convergence defined against the class of test functions given by tensor products of bounded measurable functions on the individual coordinate-spaces.  It is standard that this topology on the convex set of couplings is compact: see, for instance, Theorem 6.2 of Glasner~\cite{Gla03}.

However, we will actually deduce Theorem~\ref{thm:main} from a stronger kind of convergence.  For any joining $\l$ and any fixed choice of $f_i \in L^\infty(\mu_i)$ for $1 \leq i \leq k$, the map
\[f_0\mapsto \int_{X_0\times X_1\times \cdots \times X_k} f_0\otimes f_1\otimes \cdots \otimes f_k\,\d\l\]
defines a bounded linear functional on $L^2(\mu_0)$, and hence by the self-duality of Hilbert space it specifies a function
\[M^\l(f_1,\ldots,f_k) \in L^2(\mu_0)\]
(an alternative, more concrete description of $M^\l$ can be found in Section~\ref{sec:background} below).  The joining convergence asserted by Theorem~\ref{thm:main} is equivalent to the weak convergence in $L^2(\mu_0)$ of the averages
\[A^\l_T(f_1,\ldots,f_k) := \barint_0^T M^\l\big(f_1\circ u_1^{\phi_1(t)},\ldots,f_k\circ u_k^{\phi_k(t)}\big)\,\d t,\]
but in fact the methods we call on below (particularly the van der Corput estimate, Lemma~\ref{lem:vdC}) naturally give more:

\begin{thm}\label{thm:main2}
In the setting of Theorem~\ref{thm:main}, the averages $A^\l_T(f_1,\ldots,f_k)$ converge in norm in $L^2(\mu_0)$ as $T\to\infty$ for any functions $f_i \in L^\infty(\mu_i)$, $1 \leq i \leq k$.
\end{thm}

Of course, this does not immediately imply the remainder of Theorem~\ref{thm:main} concerning the extra symmetries of the limit joining.  That will require some additional argument.

The problem of pointwise convergence of the averages $A^\l_T$ remains open, and the methods of the present paper probably say very little about it.  One related special case (for certain discrete-time averages) has been established by Bourgain in~\cite{Bou90}, but I know of no more recent extensions of his work.

\subsection*{Origin and relation to other works}

Theorem~\ref{thm:main} has its origin in the study of multiple recurrence.  Furstenberg's original Multiple Recurrence Theorem~\cite{Fur77} asserts that for a single probability-preserving transformation $T\actson (X,\S,\mu)$, if $A \in \S$ has $\mu(A) > 0$ then also
\begin{eqnarray}\label{eq:multirec}
\liminf_{N\to\infty}\frac{1}{N}\sum_{n=1}^N \mu(A\cap T^{-n}A\cap \cdots \cap T^{-(k-1)n}A) > 0\quad\quad\forall k\geq 1.
\end{eqnarray}
In particular, there must be a time $n \geq 1$ at which
\[\mu(A\cap T^{-n}A\cap \cdots \cap T^{-(k-1)n}A) > 0:\]
this is `$k$-fold multiple recurrence' for $A$.

Furstenberg studied this phenomenon in order to give a new proof of a deep theorem of Szemer\'edi in additive cominatorics~\cite{Sze75}, which can be deduced quite easily from the Multiple Recurrence Theorem.  Following Furstenberg's original paper, many other works have either proved analogous multiple recurrence assertions in more general settings or analysed the `multiple' ergodic averages of the kind appearing in~(\ref{eq:multirec}), in particular to determine whether they converge.  We will not attempt to give complete references here, but refer the reader to~\cite{Aus--thesis}, to the paper~\cite{HosKra05} of Host and Kra and to Chapters 10 and 11 of Tao and Vu's book~\cite{TaoVu06} for more details.

Many of these convergence questions can be phrased in terms of convergence of joinings, much in the spirit of Theorem~\ref{thm:main}.  In Furstenberg's original setting, if we let $\mu^\Delta$ be the copy of $\mu$ supported on the diagonal in $X^k$, then the above averages may be re-written as
\[\int_{X^k} 1_A\otimes 1_A\otimes \cdots\otimes 1_A\,\d\mu_N,\]
where
\[\mu_N := \frac{1}{N}\sum_{n=1}^N (\id_X\times T\times \cdots\times T^{k-1})^n_\ast \mu^\Delta,\]
so in fact the convergence of these scalar averages is almost precisely the assertion that the orbit of the joining $\mu^\Delta$ under the off-diagonal $\id_X\times T\times\cdots\times T^{k-1}$ is equidistributed relative to some limit joining.  Convergence here follows from work of Host and Kra~\cite{HosKra05} (see also Ziegler~\cite{Zie07}), and it is worth noting that in this situation the additional invariance of the limit joining under $\id_X\times T\times \cdots\times T^{k-1}$ is obvious from the definition of the $\mu_N$ and the F\o lner property of the intervals $\{1,2,\ldots,N\} \subset \bbZ$.

On the other hand, that additional invariance can be put at the heart of an alternative proof of convergence, which also applies to the more general question of the convergence of the averaged joinings
\[\frac{1}{N}\sum_{n=1}^N (\id_X\times T_1\times \cdots\times T_k)^n_\ast \mu^\Delta\]
for a commuting tuple of transformations $T_1,T_2,\ldots,T_k\actson (X,\S,\mu)$: see~\cite{Aus--nonconv,Aus--thesis} (and compare with Tao~\cite{Tao08(nonconv)}, where the first proof of convergence for this higher-rank setting was given using very different methods).  This more general setting still exhibits a multiple recurrence phenomenon with striking combinatorial consequences, as shown much earlier by Furstenberg and Katznelson~\cite{FurKat78}. Another aspect of the study of the limit of the above joinings is that a sufficiently detailed understanding of its structure can be used to give an alternative proof of their theorem~\cite{Aus--newmultiSzem}.

Having come this far, it is natural to ask after the behaviour of these averaged joinings if $T_1$, $T_2$, \ldots, $T_k$ do not commute, but generate some more complicated discrete group.  In particular, if they generate a nilpotent group, then Leibman has shown that multiple recurrence phenomena still occur~\cite{Lei98} using an extension of Furstenberg and Katznelson's arguments, but that approach does not prove that the associated functional averages converge in $L^2(\mu)$.  The question of convergence seems to be closely related to whether the averages
\[\frac{1}{N}\sum_{n=1}^N (\id_X\times T^{p_1(n)}\times \cdots\times T_k^{p_k(n)})_\ast \mu^\Delta\]
converge for a $\bbZ^d$-action $T$ and polynomials $p_i:\bbZ\to\bbZ^d$, at least insofar as some of the standard methods in this area (particularly the van der Corput estimate) run into very similar difficulties in the contexts of these two problems.

These more general convergence questions were posed by Bergelson as Question 9 in~\cite{Ber96}, having previously been popularized by Furstenberg.  Several special cases were established in~\cite{FurWei96,BerLei02,HosKra05poly,Lei05(poly),Aus--lindeppleasant2,ChuFraHos09}.  On the other hand, the paper~\cite{BerLei04} contains an example in which $k=2$, $\langle T_1,T_2\rangle$ is a two-step solvable group, and convergence fails.

Shortly before the present paper was submitted, Miguel Walsh offered in~\cite{Wal11} a proof of convergence for general nilpotent groups and tuples of polynomial maps, so answering the question of Furstenberg and Bergelson in full generality.  His proof is most akin to Tao's convergence proof in~\cite{Tao08(nonconv)}, but clearly involves some non-trivial new ideas as well.  It is quite different from the very `structural' approach taken by most ergodic theoretic papers, such as the present one.  It seems likely that Walsh's approach can be adapted to prove convergence in our setting (Theorem~\ref{thm:main2}), but it gives much less information on the structure of the resulting factors and joinings (as, for example, in the rest of Theorem~\ref{thm:bigmain}).

Our Theorem~\ref{thm:main2} establishes the analog of the conjecture of Furstenberg and Bergelson (involving both nilpotent groups and polynomial maps) for continuous-time flows.  In Subsection~\ref{subs:compare-discrete} we will offer some discussion of the additional difficulties presented by an adaptation of our approach to the discrete-time setting.  It would still be of interest to find a successful such adaptation, since it would presumably require uncovering a more detailed description of the relevant factors and joinings, and so would comprise a substantial complement to the approach via Walsh's methods.

We should note also that the case $G = \bbR^d$ in Theorems~\ref{thm:main} and~\ref{thm:main2} was recently established in~\cite{Aus--ctspolyaves}.  However, the methods below diverge quite sharply from that previous paper. That work relied crucially on making a time change $t \mapsto t^\a$ in the integral averages under study for some small $\a > 0$, in order to convert averages along polynomial orbits into averages along orbits given by a linear map perturbed by some terms that grow at sublinear rates in $t$.  That trick leads to a substantial simplification of the necessary induction on families of polynomials (in that paper Bergelson's PET induction is not needed, since something more direct suffices, whereas this induction scheme will appear in the present paper shortly), and so cuts out various other parts of the argument that we use below.  However, I do not know how to implement this time-change trick for maps into general nilpotent groups, essentially because various commutators that appear during the proof can give rise to high-degree terms which disrupt the choice of any particular $\a$ used to make the leading-order terms linear.  It is also my feeling that the argument given below reveals rather more about the relevant structures within probability-preserving $G$-actions that are responsible for the asymptotic behaviour of the averages in Theorem~\ref{thm:main}.

Although it emerges from the study of multiple recurrence, Theorem~\ref{thm:main} fits neatly into the general program of equidistribution.  Equidistribution phenomena for sequences in compact spaces, and especially sequences arising from dynamical systems, have been popular subjects of analysis for most of the twentieth century: see, for instance, the classic text~\cite{KuiNie74}.  Theorem~\ref{thm:main} can be seen as a close analog of more classical results concerning special classes of compact topological systems: in place of the orbit of an individual point or distinguished subset, we study the orbit of an initially-given joining, and correspondingly vague convergence of measures (that is, tested against continuous functions on a compact space) is replaced by convergence in the coupling topoology.

Of course, equidistribution theorems for topological systems always rely very crucially on the special structure of the system under study.  Among arbitrary actions on compact spaces there are plentiful examples for which the set of invariant probabilities is very large and unstructured, and which have many points that do not equidistribute.  It is interesting that once a tuple of systems $(X_i,\S_i,\mu_i,u_i)$ with invariant probabilities has been fixed, their joinings exhibit the behaviour of Theorem~\ref{thm:main} without any extra assumptions on those individual systems.  Instead, the necessary provisions are that we start with the orbit of some joining, rather than of a single point, and then prove equidistribution in the sense of the coupling topology. 

Among the most profound results giving equidistribution for concrete systems are those concerning the orbits of unipotent flows on homogeneous spaces.  In this setting the heart of such an analysis is typically a classification of all invariant probability measures on a system, which then restricts the possible vague limits one can obtain from the empirical measures along an orbit of the system so that, ideally, one can prove that the empirical measures have only one possible limit (and so are equidistributed).

To some extent the approach to Theorem~\ref{thm:main} parallels that strategy, in that the additional invariances of the limit joinings are an important tool in the proof, and our arguments do imply some further results on the possible structure of the limit joinings (see the second remark following Proposition~\ref{prop:vdC-appn}).

The full strength of measure classification for probabilities on homogeneous spaces that are invariant and ergodic under the action of a subgroup generated by unipotent elements was finally proved by Ratner in~\cite{Rat91-a,Rat91-b}, building on several important earlier works of herself and others.  The monograph~\cite{Mor05} gives a thorough account of this story.  Following Ratner's work, Shah proved in~\cite{Sha94} some equidistribution results for trajectories of points in homogeneous spaces under flows given by regular algebraic maps into the acting group.  That notion of `polynomial' encompasses ours in many cases, and so his work offers a further point of contact between the two settings.

However, the details of the arguments used below are rather far from those developed by Ratner and her co-workers.  For instance, in Shah's paper, he first shows that any vague limit measure for the trajectory of a point under one of his regular algebraic maps must have some invariance under a nontrivial unipotent subgroup. In light of this he can restrict his attention to the possible limit measures that are permitted by Ratner's Measure Classification Theorem, whereupon the extra analysis needed can proceed. By contrast, it is essential in our work that we allow general polynomial maps into $G$ throughout, since our induction would not remain among homomorphisms even if we started there.  It would be interesting to know whether an alternative approach to Theorem~\ref{thm:main} can be found which is more in line with those works on homogeneous space dynamics.

\subsection*{First outline of the proof}

Theorems~\ref{thm:main} and~\ref{thm:main2} will be proved by induction on the tuple of polynomial maps $(\phi_1,\phi_2,\ldots,\phi_k)$.  The ordering on polynomials that organizes this induction is (a variant of) Bergelson's PET ordering from~\cite{Ber87}, which has become a mainstay of the study of multiple averages involving nilpotent groups or polynomial maps.

To a large extent, the new innovation below is the formulation of an assertion that includes Theorem~\ref{thm:main} and can be closed on itself in this induction.  The delicacy of this formulation is largely attributable to the van der Corput estimate (Lemma~\ref{lem:vdC}), which relates the averages involving a given tuple of polynomial maps to another tuple that precedes it in the PET ordering.  In the first place, it is this that forces us to prove Theorem~\ref{thm:main2} alongside Theorem~\ref{thm:main}, but it will also required other features in our inductive hypothesis.

An application of this lemma converts an assertion about a tuple of polynomial maps
\[t\mapsto \phi_i(t)\]
into another about the `differenced' maps
\begin{eqnarray}\label{eq:diffs}
(t,s)\mapsto \phi_i(t+s)\phi_i(t)^{-1}
\end{eqnarray}
(or more complicated relatives of these: see Section~\ref{sec:poly}). Regarded as functions of $t$ alone, these precede the tuple $(\phi_1,\ldots,\phi_k)$ in the PET ordering for any fixed $s$.  In many applications of PET induction one simply forms these derived maps, then fixes a value of $s$ and applies an inductive hypothesis to the restrictions of these new maps to $\bbR\times \{s\}$.  Unfortunately, in our setting there can be some values of $s$ for which the behaviour of these restrictions is not as `good' as our argument needs.  To overcome this we must retain the picture of the new maps in~(\ref{eq:diffs}) as being polynomial on the whole of $\bbR\times \bbR$.  As a consequence of this polynomial structure and certain general results about actions of nilpotent Lie groups (see Section~\ref{sec:nil-actions}), one finds that these averages behave `asymptotically the same' for all but a small set of exceptional values of $s$.  This turns out to be a crucial improvement over the possible worst-case behaviour over $s$.  Since repeated appeals to the van der Corput estimate lead to a proliferation of these differencing parameters $s$, we must actually formulate a theorem which allows for polynomial maps $\bbR\times \bbR^r\to G$, where we average over the first coordinate in $\bbR\times \bbR^r$ and the theorem promises some additional good behaviour for generic values of the remaining $r$ coordinates.

The right notion of genericity to make this precise is provided by Baire's definition of category, but transplanted into the Zariski topology of $\bbR^n$ (which is not Hausdorff and so not quite in the usual mould for applications of Baire category).  The required notion of `Zariski genericity' will be defined in Section~\ref{sec:Zariski}, and will be found to relate very well to other standard notions of `smallness' for subsets of $\bbR^n$.

In terms of this definition, the complete statement that will be proved by PET induction is as follows.

\begin{thm}\label{thm:bigmain}
Suppose that $(X_i,\S_i,\mu_i,u_i)$, $0 \leq i \leq k$, and $\l$ are as above and that $\phi_i:\bbR\times \bbR^r\to G$, $1 \leq i \leq k$, are polynomial maps satisfying $\phi_i(0,\cdot) \equiv e$.  Let $M^\l$ be constructed from $\l$ as previously, let
\[A^\l_T(f_1,\ldots,f_k) := \barint_0^T M^\l\big(f_1\circ u_1^{\phi_1(t,h)},\ldots,f_k\circ u_k^{\phi_k(t,h)}\big)\,\d t\]
(so $A^\l_T$ implicitly depends on $h$), and let
\[\vec{\phi}:= (e,\phi_1,\phi_2,\ldots,\phi_k):\bbR\times\bbR^r\to G^{k+1}.\]
Then
\begin{enumerate}
\item for any $h \in \bbR^r$ and any $f_i \in L^\infty(\mu_i)$, $1 \leq i \leq k$, the functional averages $A^\l_T(f_1,\ldots,f_k)$
converge in $L^2(\mu_0)$ as $T \to\infty$,
\item for any $h \in \bbR^r$ the averaged joinings
\[\barint_0^T (\id_{X_0}\times u_1^{\phi_1(t,h)}\times u_2^{\phi_2(t,h)}\times \cdots \times u_k^{\phi_k(t,h)})_\ast\l\,\d t\]
converge as $T \to \infty$ to some limit joining $\l^h$ which is invariant under
\[\langle G^{\Delta (k+1)}\cup \img \vec{\phi}(\cdot,h)\rangle,\]
and
\item the map $h\mapsto \l^h$ is Zariski generically constant on $E$, and the generic value it takes is a joining invariant under
\[\langle G^{\Delta (k+1)}\cup \img\vec{\phi}\rangle.\]
\end{enumerate}
\end{thm}

This clearly implies both of the previous theorems.  The rest of the paper is directed towards the proof of Theorem~\ref{thm:bigmain}.

\subsection*{Overview of the paper}

Sections~\ref{sec:background} through~\ref{sec:idem} establish certain background results that we will need for the main proofs, concerning general properties of group actions and representations; polynomial maps and genericity in the Zariski topology; finer results about actions of nilpotent Lie groups; and the technology of `idempotent' classes of probability-preserving systems.  Once all this is at our disposal, the proof of Theorem~\ref{thm:bigmain} is completed in Sections~\ref{sec:k=2},~\ref{sec:char-factor} and~\ref{sec:general-k}.  Finally Section~\ref{sec:further-ques} contains a discussion of various further questions related to those of this paper.

\section{Background on group actions}\label{sec:background}

If $G$ is a locally compact second countable (`l.c.s.c.') group, then a \textbf{$G$-system} is a quadruple $(X,\S,\mu,u)$ in which $(X,\S,\mu)$ is a standard Borel probability space and $g\mapsto u^g$ is a jointly measurable, $\mu$-preserving left action of $G$ on $(X,\S)$.  Sometimes this situation will alternatively be denoted by $u:G\actson (X,\S,\mu)$, and sometimes a whole system will be denoted by a boldface letter such as $\bfX$.

Relatedly, a \textbf{$G$-representation} is a strongly continuous orthogonal representation $\pi$ of $G$ on a separable real Hilbert space $\frH$.  (It would be more conventional to work with complex Hilbert spaces and unitary representations, but choosing the real setting avoids the need to keep track of several complex conjugations later.) This situation will often be denoted by $\pi:G\actson\frH$.  Given a $G$-system $(X,\S,\mu,u)$, the associated \textbf{Koopman representation} $u^\ast:G\actson L^2(\mu)$ is defined by
\[u(g)^\ast f := f\circ u^{g^{-1}},\]
where this convention concerning inverses ensures that both $u$ and $u^\ast$ are left actions. Here and throughout the paper the notation $L^p$, $1 \leq p \leq \infty$, is used for real Lebesgue spaces. It is classical that the joint measurability of $u$ implies the strong continuity of $u^\ast$ (see, for instance, Lemma 5.28 of Varadarajan~\cite{Varadara85}), so the Koopman representation is a $G$-representation in the present sense.

Given a $G$-system and a closed subgroup $H\leq G$, one may construct the $\s$-subalgebra
\[\S^H := \{A \in \S:\ \mu(u^hA \triangle A) = 0\ \forall h\in H\}.\]
If $H$ is normal in $G$ then this is globally $G$-invariant, and hence defines a factor of the original system which we call the \textbf{$H$-partially invariant} factor.  For some quite special technical reasons we will need only the case of normal $H$ in this paper: see Corollary~\ref{cor:normal-clos} below.

Similarly, for a $G$-representation $\pi$ we let
\[\rm{Fix}(\pi(H)):= \{v \in \frH:\ \pi(h)v = v\ \forall h \in H\} \leq \frH;\]
for Koopman representations it is easily seen that
\[\rm{Fix}(u^\ast(H)) = L^2(\mu|_{\S^H}).\]

Sometimes it is necessary to compare actions of different groups.  If $q:H\to G$ is a continuous homomorphism of l.c.s.c. groups and $\bfX = (X,\S,\mu,u)$ is a $G$-system, then we may define an $H$-system on the same probability space by letting $h$ act by $u^{q(h)}$. We denote this system by $\bfX^{q(\cdot)} = (X,\S,\mu,u^{q(\cdot)})$.  A similar construction is clearly possible for representations.

We will also need certain standard calculations involving couplings and joinings.  Suppose that $\l$ is a coupling of $\mu_0$, $\mu_1$, \ldots, $\mu_k$ (without any assumption about group actions).  We may regard it instead as a coupling of $(X_0,\S_0,\mu_0)$ with
\[(X_1\times \cdots\times X_k,\S_1\otimes\cdots\otimes \S_k,\l')\]
where $\l'$ is the marginal of $\l$ on the last $k$ coordinates.  Now $\l$ can be disintegrated over the first coordinate to obtain a probability kernel
\[\L:X_0 \to \rm{Pr}(X_1\times \cdots\times X_k,\S_1\otimes\cdots\otimes \S_k)\]
so that
\[\l = \int_{X_0}\delta_{x_0}\otimes \L(x_0,\,\cdot\,)\,\mu_0(\d x_0);\]
and this, in turn, defines a multilinear map
\[M^\l:L^\infty(\mu_1)\times \cdots\times L^\infty(\mu_k)\to L^\infty(\mu_0)\]
according to
\[M^\l(f_1,\ldots,f_k)(x_0) := \int_{X_1\times \cdots\times X_k}f_1 \otimes f_2\otimes \cdots \otimes f_k\,\d\L(x_0,\,\cdot\,).\]
Clearly one has
\[\int_{X_0}f_0\cdot M^\l(f_1,\ldots,f_k)\,\d\mu_0 = \int_{X_0\times X_1\times \cdots\times X_k}f_0\otimes f_1 \otimes \cdots \otimes f_k\,\d\l,\]
so this agrees with the definition of $M^\l$ by duality given in the Introduction.

The following is now a routine re-formulation of the definition of a relatively independent product, and the proof is omitted; see, for instance, the third of Examples 6.3 in Glasner~\cite{Gla03}.

\begin{lem}\label{lem:l_1l}
Let $\L:X_0\to \Pr(X_1\times \cdots \times X_k)$ be as above and define the relative product measure $\l\otimes_0\l$ on $X_1^2\times \cdots \times X_k^2$ by
\[\l\otimes_0\l = \int_{X_0}\L(x_0,\cdot)\otimes \L(x_0,\cdot)\,\mu_0(\d x_0).\]
Then for any $f_i,g_i \in L^\infty(\mu_i)$, $1 \leq i \leq k$, one has
\begin{multline*}
\int_{X_0}M^\l(f_1,f_2,\ldots,f_k)\cdot M^\l(g_1,g_2,\ldots,g_k)\,\d\mu_0\\ = \int_{X_1^2\times\cdots\times X_k^2}f_1\otimes g_1\otimes f_2\otimes g_2\otimes \cdots \otimes f_k\otimes g_k\,\d (\l\otimes_0\l).
\end{multline*}
\qed
\end{lem}

\section{Real polynomials and Zariski residual sets}\label{sec:Zariski}

The third part of Theorem~\ref{thm:bigmain} involves the notion of Zariski genericity.  Recall that on $\bbR^n$ (or any other real algebraic variety) the \textbf{Zariski topology} is the topology whose closed sets are the subvarieties.  Although the failure of $\bbR$ to be algebraically closed gives rise to certain novel behaviour not seen in more classical algebraic geometry (especially under projection maps), in this paper we will not meet any of the situations in which this matters.  The basic notions of the theory can be found in many books that use algebraic groups, such as in Subsection D.1 of Starkov~\cite{Sta00}.  The additional idea we need from that arena is the following.

\begin{dfn}[Zariski meagre and residual sets]
A subset $W \subseteq \bbR^n$ is \textbf{Zariski meagre} if it can be covered by a countable family of proper subvarieties of $\bbR^n$.  A subset of $\bbR^n$ is \textbf{Zariski residual} if its complement is Zariski meagre.  A property that depends on a parameter $h \in \bbR^n$ is \textbf{Zariski generic} if it obtains on a Zariski residual set of $h$.
\end{dfn}

Since proper subvarieties are always closed and nowhere dense in the Euclidean topology, Zariski residual sets are residual in the Euclidean topology.  They are therefore `large' in the sense of the Baire Category Theorem and its consequences, but in a much more structured way than an arbitrary Euclidean-residual subset.  In particular, they exhibit the following simple behaviour under slicing:

\begin{lem}\label{lem:Zar-res-on-subspace}
If $E \subseteq \bbR^n$ is Zariski meagre and $V \subseteq \bbR^n$ is any affine subspace then either $E \supseteq V$ or $E \cap V$ is Zariski meagre in $V$.  In the space of translates $\bbR^n/V$, the subset of translates for which the former holds is Zariski meagre.
\end{lem}

\textbf{Proof}\quad This is simply a consequence of the corresponding property of Zariski closed sets. \qed

Zariski meagre sets are also small in a natural measure-theoretic sense.

\begin{lem}
A Zariski meagre subset $E \subseteq \bbR^n$ has Hausdorff dimension at most $n - 1$.
\end{lem}

\textbf{Proof}\quad Clearly it suffices to show that a single proper algebraic subvariety $V \subseteq \bbR^n$ has Hausdorff dimension at most $n - 1$, and moreover that this holds when $V = \{f = 0\}$ for some nonzero polynomial $f:\bbR^n\to \bbR$ (because any proper $V$ can be contained in such a zero-set).

This follows by induction on degree.  If $f$ is linear then it is immediate, so suppose $\deg f \geq 2$. Then on the one hand the nonsingular locus $\{f = 0\}\cap \{\nabla f \neq 0\}$ can be covered with countably many open sets on which $\{f = 0\}\cap \{\nabla f\neq 0\}$ locally agrees with a smooth $(n-1)$-dimensional submanifold of $\bbR^n$, and hence has Hausdorff dimension $n-1$. On the other hand, the remaining set $\{f = 0\}\cap \{\nabla f = 0\}$ is contained in the set $\{\ell(\nabla f) = 0\}$ for any choice of $\ell \in (\bbR^n)^\ast\setminus \{0\}$, which is an algebraic variety generated by a polynomial of degree at most $\deg f - 1$ and so has Hausdorff dimension at most $n-1$ by the inductive hypothesis. \qed

\section{Polynomial maps into nilpotent Lie groups}\label{sec:poly}

Henceforth $G$ will denote a connected and simply connected nilpotent Lie group, $\frg$ its Lie algebra,
\[G = G^1 \unrhd G^2 \unrhd \cdots \unrhd G^s\unrhd (e)\]
its ascending central series, and
\[\frg = \frg^1 \unrhd \frg^2\unrhd \cdots \unrhd \frg^s\unrhd (0)\]
the corresponding ascending series of $\frg$.

In the following we will need certain standard facts about such groups, in particular that the exponential map $\exp:\frg\to G$ is an analytic diffeomorphism and that any Lie subalgebra $\frh \leq \frg$ exponentiates to a closed Lie subgroup of $G$, which is normal if and only if $\frh$ was an ideal.  (Note that both of these require the assumption that $G$ is simply connected as well as connected.)  These can be found as Theorem 1.2.1 and Corollary 1.2.2 in Corwin and Greenleaf~\cite{CorGre90}, which provides a good general reference for the study of these groups.

\subsection{Polynomial maps}

\begin{dfn}[Polynomial map]\label{dfn:poly}
A map $\phi:G'\to G$ between nilpotent Lie groups is \textbf{polynomial} if there is some $d \geq 1$ such that
\[\nabla_{h_1}\nabla_{h_2}\cdots\nabla_{h_d}\phi \equiv e \quad\quad \forall h_1,h_2,\ldots,h_d \in G',\]
where $\nabla_h \phi(g) := \phi(gh^{-1})\phi(g)$.
\end{dfn}

This definition has come to prominence in the study of multiple recurrence phenomena since Leibman's work generalizing the Furstenberg-Katznelson Multiple Recurrence Theorem to tuples of transformations generating a nilpotent group~\cite{Lei98}.  For maps into a module $M$ over a ring $R$ (such as an Abelian group, which is a module over $\bbZ$), degree-$d$ polynomial maps have been studied much more classically as an ideal of functions $G\to M$ annihilated under convolution by the $d^{\rm{th}}$ power of the augmentation ideal of $R[G]$: see, for instance, Passi~\cite{Pas68,Pas79}.

In this work we will need the above definition only for $G' = \bbR^n$.  If in addition $G = \bbR^m$, then it is a simple exercise to show that a map $\phi$ is polynomial according to the above if and only if it may be expressed as an $m$-tuple of polynomials in $n$ variables. For general nilpotent targets $G$ a more concrete view of polynomial maps is still available by the following standard proposition and corollary (for the former see, for instance, Proposition 1.2.7 in Corwin and Greenleaf~\cite{CorGre90}).

\begin{prop}
If $G$ is an $s$-step connected and simply connected nilpotent Lie group, then $\exp:\frg \to G$ is a diffeomorphism, and pulled back through $\exp$ the operations of multiplication and inversion become polynomial maps $\frg\times\frg\to \frg$ and $\frg\to \frg$ of degree bounded only in terms of $s$. \qed
\end{prop}

\begin{cor}\label{cor:exppoly-is-poly}
A map $\phi:\bbR^n\to G$ is polynomial if and only if it is of the form $\exp\circ \Phi$ for some polynomial $\Phi:\bbR^n\to \frg$.
\end{cor}

\textbf{Proof}\quad This follows by induction on the nilpotency class of $G$. On the one hand, if $\Phi:\bbR^n\to \frg$ is a polynomial, then after $(\deg\Phi)$-many applications of the differencing operator $\nabla_\bullet$ the exponentiated map $\exp\circ \Phi$ may not vanish identically, but at least its projection to $G/G^2$ vanishes because this is isomorphic to the projection of $\Phi$ to $\frg/\frg^2$.  Thus finitely many differencing operations yield a polynomial map into $\frg^2$, and now repeating this argument $s$ times shows that the differences of $\exp\circ\Phi$ do eventually vanish.

On the other hand, if $\phi:\bbR^n \to G$ is a polynomial map, then the same is true of $\phi G_2:\bbR^n\to G/G^2\cong \bbR^{\dim G- \dim G^2}$.  This, in turn, is simply isomorphic to $(\exp^{-1}\circ\phi) + \frg^2:\bbR^n\to \frg/\frg^2$, so this latter is a polynomial.  By choosing lifts of its coefficients under the projection $\frg\to\frg/\frg^2$, we obtain a polynomial $\Phi_1 :\bbR^n\to\frg$ such that $\exp\circ(\exp^{-1}\circ \phi - \Phi_1)$ takes values in $G^2$, and it is clearly still a polynomial map there using the argument of the previous paragraph.  Now the inductive hypothesis applied to $G^2$ gives another polynomial $\Phi_2:\bbR^n\to \frg^2$ such that $\exp\circ(\exp^{-1}\circ \phi - \Phi_1) = \exp\circ \Phi_2$, and re-arranging this completes the proof. \qed

By pulling back to the Lie algebra and arguing there, the above proposition and corollary have the following further consequence, which will be useful in the sequel.

\begin{cor}
If $\phi,\psi:\bbR^n\to G$ are polynomial maps, then so are the pointwise product $x\mapsto \phi(x)\psi(x)$ and the pointwise inverse $x\mapsto \phi(x)^{-1}$. \qed
\end{cor}

\subsection{Families of maps and the PET ordering}

Our attention now turns to finite tuples
\[\F = (\phi_1,\phi_2,\ldots,\phi_k)\]
of polynomial maps $\bbR\times\bbR^r\to G$.

In what follows it is extremely important that we consider the domain of these maps to be split as $\bbR\times \bbR^r$.  Although this is not really different from $\bbR^{r+1}$, the heart of the main induction below rests on comparing the degrees of different polynomial maps into $G$ \emph{in the first coordinate only}.  Therefore we will henceforth restrict attention to maps defined on products of $\bbR$ with other real vector spaces, and will always regard the second coordinate as an auxiliary parameter.

\begin{dfn}[Internal class; leading degree; leading term]
For a polynomial map $\phi:\bbR \times \bbR^r\to G$ with $\phi(0,\cdot) \equiv e$, its \textbf{internal class} is the greatest $c$ such that $G^c \supseteq \img \phi$.  It is denoted $\rm{cl}\,\phi$.

Given this, the projection
\[\phi G^{c+1}:(t,h)\mapsto \phi(t,h)G^{c+1}:\bbR\times \bbR^r \to G^c/G^{c+1} \cong \bbR^{\dim G^c - \dim G^{c+1}}\]
is a Euclidean-valued polynomial map.  The \textbf{leading degree} $\ldeg \phi$ of $\phi$ is the degree of $\phi G^{c+1}$ in the variable $t$, and the \textbf{leading term} of $\phi$ is the term in $\phi G^{c+1}$ of the form $t^{\ldeg \phi}\psi(h)$ for some polynomial map $\psi:\bbR^r\to G^c/G^{c+1}$.
\end{dfn}

\begin{dfn}[Leading-term equivalence]
Two polynomial maps $\phi,\psi:\bbR\times \bbR^r\to G$ are \textbf{leading-term equivalent}, denoted $\phi \sim_{\LT} \psi$, if $\rm{cl}\,\phi = \rm{cl}\,\psi$ and $\phi$ and $\psi$ have the same leading term (hence certainly the same leading degree).
\end{dfn}

Several further definitions are needed in order to explain the PET ordering that will steer the inductive proof of Theorem~\ref{thm:bigmain}.  The next roughly follows Leibman~\cite{Lei98}.

\begin{dfn}[Weight]
The \textbf{weight} of a polynomial $\phi:\bbR\times \bbR^r\to G$ is the pair $\wt\phi := (\rm{cl}\phi,\ldeg\phi)$. The set $\Wt$ of possible weights $(c,d)$ is ordered lexicographically: pairs $(c,d),(c',d') \in \Wt$ satisfy $(c,d) \prec (c',d')$ if
\begin{itemize}
\item either $c > c'$,
\item or $c = c'$ and $d < d'$.
\end{itemize}
Since clearly $\phi\sim_\LT\psi$ implies $\wt\phi = \wt\psi$, we may also define the \textbf{weight} of an $\sim_\LT$-equivalence class as the weight of any of its members.
\end{dfn}

This is a well-ordering on $\Wt$, and it now gives rise to a partial ordering on polynomial maps.

\begin{dfn}[PET ordering on polynomials]
Given two polynomial maps $\phi, \psi:\bbR\times \bbR^r\to G$, the first \textbf{precedes} the second in the \textbf{PET ordering}, denoted $\phi\prec_\PET\psi$, if $\wt\phi\prec \wt\psi$.
\end{dfn}

\textbf{Remark}\quad Our $\prec_\PET$ is not quite the same as the PET ordering used in much of the earlier literature for polynomial maps into nilpotent groups.  Those required a comparison between polynomials in terms of the individual members of some Mal'cev basis of $G$; see, for instance, Section 3 in~\cite{Lei98}.  Our ordering is actually a little weaker (in the sense that $\prec_\PET \subsetneqq \prec_\PET^{\rm{previous}}$ as relations), because we compare our polynomials on the whole Euclidean subquotients of $G$ arising from the ascending central series, and so in our ordering the assertion that two polynomials have the same leading term is stronger.  However, when we later use the PET induction via the van der Corput lemma it will be clear that we are still moving strictly downwards among our families of polynomials, so that the induction proceeds correctly. \fin

The PET ordering on polynomials will play a r\^ole in the proof of the special case $k=2$ of Theorem~\ref{thm:bigmain}, but the general case will require an extension of it to an ordering of tuples of polynomials.

\begin{dfn}
Suppose that $f,g:\Wt\to\bbN$ are maps which each take nonzero values at only finitely many weights.  Then $f$ \textbf{precedes} $g$, denoted $f \prec g$, if there is some $(c,d) \in \Wt$ such that
\begin{itemize}
\item $f(c',d') = g(c',d')$ whenever $(c',d') \succ (c,d)$, and
\item $f(c,d) < g(c,d)$.
\end{itemize}
\end{dfn}

\begin{dfn}[PET ordering for tuples of polynomials]
If $\F = (\phi_1,\phi_2,\ldots,\phi_k)$ is a tuple of polynomial maps then its \textbf{weight assignment} is the function $\Wt\F:\Wt \to \bbN$ which to each $(c,d) \in \Wt$ assigns the number of $\sim_{\rm{LT}}$-equivalence classes of maps in $\F$ that have weight $(c,d)$.

Suppose now that $\F = (\phi_1,\phi_2,\ldots,\phi_k)$ and $\cal{G} = (\psi_1,\psi_2,\ldots,\psi_\ell)$ are families of polynomial maps $\bbR\times \bbR^r\to G$.  Then $\F$ \textbf{precedes} $\cal{G}$, denoted $\F \prec_\PET \cal{G}$, if
\begin{itemize}
\item either $\Wt\F\prec\Wt\cal{G}$,
\item or $\Wt\F = \Wt\cal{G}$, and the sets of $\sim_\LT$-equivalence classes $\F/\!\!\sim_\LT$ and $\cal{G}/\!\!\sim_\LT$ can be matched in such a way that (i) their weights match, (ii) every class of $\F$ has cardinality no larger than its corresponding class in $\cal{G}$, and (iii) in at least one instance it is strictly smaller.
\end{itemize}
\end{dfn}

As in most proofs that use the PET ordering, it is needed for a particular pair of families of maps, one derived from the other according to the following definitions.

\begin{dfn}[Pivot]
If $\F = (\phi_1,\phi_2,\ldots,\phi_k)$ is a tuple of polynomial maps $\bbR\times \bbR^r\to G$ then a \textbf{pivot} for $\F$ is a PET-minimal member $\phi \in \F$.
\end{dfn}

\begin{dfn}[Derived family]\label{dfn:derived}
Suppose that $\F = (\phi_1,\phi_2,\ldots,\phi_k)$ is a tuple of polynomial maps $\bbR\times \bbR^r\to G$.  Then for $i\leq k$ its \textbf{$i^{\rm{th}}$ derived family} consists of the following polynomial maps $\bbR\times (\bbR\times \bbR^r)\to G$:
\[(t,k,h)\mapsto \phi_j(t,h)\phi_i(t,h)^{-1}\quad\quad\hbox{for}\ j \in \{1,2,\ldots,k\}\setminus \{i\}\]
and
\[(t,k,h)\mapsto \phi_j(k,h)^{-1}\phi_j(t+k,h)\phi_i(t,h)^{-1}\quad\quad\hbox{for}\ j \in \{1,2,\ldots,k\}.\]
\end{dfn}

Note that the pre-multiplication by $\phi_j(k,h)^{-1}$ in the last line has the consequence that if $\phi_j(0,\cdot) \equiv e$ for every $i$, then the same is true of the derived family.

\begin{lem}\label{lem:PET-calcns}
If $\F = (\phi_1,\phi_2,\ldots,\phi_k)$ with $\phi_1$ a pivot, then its first derived family precedes it in the PET ordering.  Also, the sub-tuple $(\phi_2,\ldots,\phi_k)$ precedes $\F$ in the PET ordering.
\end{lem}

\textbf{Proof}\quad For each $j\geq 2$ consider the polynomial maps
\[\phi_j(t,h)\phi_1(t,h)^{-1}\quad\quad\hbox{and}\quad\quad \phi_j(k,h)^{-1}\phi_j(t+k,h)\phi_1(t,h)^{-1}.\]
Because $\phi_1$ is a pivot,
\begin{itemize}
\item either $\wt\phi_j \succ \wt\phi_1$,
\item or $\wt\phi_j = \wt \phi_1$ but $\phi_j \not\sim_\LT \phi_1$,
\item or $\phi_j \sim_\LT \phi_1$.
\end{itemize}
In the first case both of the new maps above still have weight equal to $\wt\phi_j$, and are actually leading-term equivalent by comparing their leading terms in $G^c/G^{c+1}$ for $c = \rm{cl}\,\phi_j$.  By the same reasoning, if $\phi_j \sim_\LT \phi_{j'}$ then all four of the resulting new maps are leading-term equivalent.

The same conclusions hold when $\wt\phi_1 = \wt\phi_j$ but $\phi_1 \not\sim_\LT \phi_j$, since in this case the leading term of either of the above maps into $G^c/G^{c+1}$ is given by the nonzero difference of the leading terms of $\phi_1$ and $\phi_j$.

Lastly, if $\phi_j \sim_\LT \phi_1$, then these leading terms do cancel, and so both of the polynomial maps written above now strictly precede $\phi_1$ in the PET ordering.

Therefore overall the equivalence classes of $\F$ and of its $1^{\rm{st}}$ derived family are in bijective weight-preserving correspondence, apart from the equivalence class of $\phi_1$, which is replaced by (possibly several) classes in the derived family of strictly lower weight. This proves the first assertion.

The second assertion is obvious, because the removal of $\phi_1$ either removes a whole $\sim_\LT$-equivalence class in case $\phi_1$ is in a singleton class, and hence reduces $\Wt\F$ in $\prec$, or leaves the $\sim_\LT$-class structure of $\F$ unchanged but reduces the cardinality of exactly one of the classes.
\qed

\section{Finer results for actions of nilpotent Lie groups}\label{sec:nil-actions}

For any inclusion $H \leq G$ of topological groups, $H^{\rm{n}}$ will denote the \textbf{topological normal closure} of $H$ in $G$: that is, the completion of the normal closure in $G$. This notation suppresses the dependence of this definition on the larger group $G$, which will always be clear from the context.  Similarly, if $G$ is a connected and simply connected Lie group with Lie algebra $\frg$ and $V \leq \frg$ is a Lie subalgebra, then $V^\rm{n}$ denotes the Lie algebra generated by $\sum_g \rm{Ad}(g)V$ (equivalently, the Lie ideal generated by $V$ in $\frg$), so that $\exp(V^\rm{n}) = (\exp V)^\rm{n}$.

The first important result we need is a consequence of the classical Mautner Phenomenon.  We will make use of the following expression of this argument as isolated by Margulis~\cite{Mar91}; it can also be found as Lemma 2.2 in Subsection 2.1 of Starkov~\cite{Sta00}.

\begin{lem}[Mautner Phenomenon]\label{lem:Mautner}
Suppose that $\pi:G\actson \frH$ is a orthogonal representation of a connected Lie group, that $H \leq G$ is a connected Lie subgroup, that $g \in G$ and that there are a sequences $g_i\in G$ and $h_i,h_i' \in H$ with $g_i \to e$ and $g_ih_ig_i^{-1}h_i'\to g$.  Then
\[\rm{Fix}(\pi(g)) \supseteq \rm{Fix}(\pi(H)).\]
\qed
\end{lem}

\begin{cor}\label{cor:normal-clos}
If $G$ is a connected and simply connected nilpotent Lie group, $H\leq G$ is a connected closed subgroup and $\pi:G\actson \frH$ is an orthogonal representation, then
\[\rm{Fix}(\pi(H)) = \rm{Fix}(\pi(H^\rm{n})).\]
Similarly, if $(X,\S,\mu,u)$ is a $G$-system then
\[\S^H = \S^{H^\rm{n}}.\]
\end{cor}

\textbf{Proof}\quad We focus on the first claim, since the second follows at once by considering the Koopman representation.

A simple calculation shows that $H^{\rm{n}} = \langle H [H,G]\rangle$, where $[H,G]$ is the subgroup generated by all commutators of elements of $H$ with elements of $G$. Let
\[G = G_1 \unrhd G_2 \unrhd \ldots \unrhd G_s \unrhd G_{s+1} = \{e\}\]
be a central series of $G$ in which each quotient $G_r/G_{r+1}$ has dimension one; for example, one may insert extra terms into the ascending central series, as in the construction of a strong Mal'cev basis.  Let $\frg_r$ be the Lie algebra of $G_r$ and $\frh$ the Lie algebra of $H$.

We will prove by downwards induction on $r$  that if $1 \leq r\leq s$ then
\[\rm{Fix}(\pi(\langle H[G_{r+1},H]\rangle)) = \rm{Fix}(\pi(\langle H[G_r,H]\rangle)).\]
When $r = s$ the left-hand side here is $\rm{Fix}(\pi(H))$, while when $r = 1$ the right-hand side is $\rm{Fix}(\pi(H^{\rm{n}}))$, so this will complete the proof.

When $r=s$ the result is clear because $G_s$ is central in $G$, so now suppose the result is known for some $r+1 \leq s$.  By replacing $H$ with $\langle H[G_{r+1},H]\rangle $, we may assume that they are equal, since another easy calculation shows that the sets
\[(H[G_{r+1},H])\cdot \big[G_{r+1},(H[G_{r+1},H])\big] \quad\quad \hbox{and} \quad\quad H[G_{r+1},H]\]
generate the same subgroup of $G$.

Let $V \in \frg_r\setminus \frg_{r+1}$, so that $\frg_r$ is the smallest Lie algebra containing both $V$ and $\frg_{r+1}$. The subgroup $\langle H[G_r,H]\rangle$ is connected, and its Lie algebra is the smallest Lie subalgebra of $\frg$ that contains both $\frh$ and $\{[V,U]:\ U \in \frh\}$.  It therefore suffices to show that any $v \in \rm{Fix}(\pi(H))$ is also fixed by $\exp([V,U])$ for any $U \in \frh$.

This can be deduced using Lemma~\ref{lem:Mautner}.  We need to show that if $U \in \frh$ then $\exp([V,U])$ is a limit of group elements of the form $g_ih_ig_i^{-1}h_i'$, as treated in that lemma.  This follows from the Baker-Campbell-Hausdorff formula, which implies for any $t > 0$ that
\[\exp(tV)\exp((1/t)U)\exp(-tV)\exp(-(1/t)U) = \exp([V,U] + \rm{O}(t))\exp(R(t)),\]
where $R(t)$ collects those multiple commutators that involve at least one copy of $V$ and at least two entries from $\frh$, which must therefore lie in
\[[\frg_{r+1},\frh]\subseteq \frh.\]
Hence
\[\exp([V,U]) = \exp(tV)\exp((1/t)U)\exp(-tV)\big(\exp(-(1/t)U)\exp(-R(t))\big),\]
so letting $t = 1/i$ gives the conditions needed by Lemma~\ref{lem:Mautner}. \qed

\begin{cor}\label{cor:rel-ind-over-common}
If $G$ is a connected and simply connected nilpotent Lie group, $H_1,H_2 \leq G$ are connected closed subgroups and $\pi:G\actson \frH$ is an orthogonal representation, then the subspaces
\[\rm{Fix}(\pi(H_1)),\ \rm{Fix}(\pi(H_2))\leq \frH\]
are relatively orthogonal over their common further subspace
\[\rm{Fix}(\pi(\langle H_1\cup H_2\rangle))\]
(meaning that
\[\rm{Fix}(\pi(H_1)) \ominus \rm{Fix}(\pi(\langle H_1\cup H_2\rangle)) \perp \rm{Fix}(\pi(H_2)) \ominus \rm{Fix}(\pi(\langle H_1\cup H_2\rangle)).\quad )\]
Similarly, if $(X,\S,\mu,u)$ is a $G$-system then $\S^{H_1}$ and $\S^{H_2}$ are relatively independent over $\S^{\langle H_1\cup H_2\rangle}$.
\end{cor}

\textbf{Proof}\quad For a Lie subgroup $H \leq G$, since
\[\rm{Fix}(\pi(H)) = \rm{Fix}(\pi(H^{\rm{n}}))\]
and $H^{\rm{n}}\unlhd G$, this subspace of $\frH$ is actually invariant under the whole action $\pi$. Therefore the orthogonal projections $P_i$ onto $\rm{Fix}(\pi(H_i))$ both commute with $\pi$.

It follows that $P_1P_2$ has image contained in $\rm{Fix}(\pi(\langle H_1\cup H_2\rangle))$. Since conversely any vector fixed by both $H_1$ and $H_2$ is also fixed by $P_1$ and $P_2$, it follows that $P_1P_2$ is an idempotent with image equal to $\rm{Fix}(\pi(\langle H_1\cup H_2\rangle))$, and the same holds for $P_2P_1$.  Hence for any vectors $u \in \frH$ and $v \in \rm{Fix}(\pi(\langle H_1\cup H_2\rangle))$ one has
\[\langle u,v\rangle = \langle u,(P_1P_2)v\rangle = \langle (P_2P_1)u,v\rangle,\]
so in fact $P_2P_1$ is the orthogonal projection onto its image, and similarly for $P_1P_2$.

Finally, if $v_i \in \rm{Fix}(\pi(H_i))$ for $i=1,2$ then this implies
\[\langle v_1,v_2\rangle = \langle P_1v_1,P_2v_2\rangle = \langle P_2P_1v_1,v_2\rangle = \langle (P_2P_1)v_1,(P_2P_1)v_2\rangle,\]
which is the desired relative orthogonality.

In the case of a $G$-system, applying the above result to the Koopman representation tells us that for any $\S^{H_i}$-measurable functions $f_i \in L^2(\mu_i)$ for $i=1,2$ we have
\[\int_X f_1f_2\,\d\mu = \int_X \sfE(f_1\,|\,\S^{\langle H_1\cup H_2\rangle})\sfE(f_2\,|\,\S^{\langle H_1\cup H_2\rangle})\,\d\mu,\]
and this is the desired relative independence. \qed

\textbf{Example}\quad The above proofs are intimately tied to the nilpotency of $G$, so it is worth including an example of a solvable Lie group $G$ and representation $\pi:G\actson \frH$ to show that this restriction is really needed.

Let $\rho:\bbR\actson\bbC$ be the rotation action defined by
\[\rho^tz := \rm{e}^{2\pi\rm{i} t}z\]
and let $G := \bbC\rtimes_\rho\bbR$.  This is a simple three-dimensional solvable Lie group; in coordinates it is $\bbC\times \bbR$ with the product
\[(u,s)\cdot (v,t) := (\rho^tu + v,s + t).\]
It may also be interpreted as a group extension of $\bbZ$ by the group $\bbC\rtimes \rm{S}^1$ of orientation-preserving isometries of $\bbC$, and this picture gives an action $\xi:G\actson \bbC$ with kernel isomorphic to $\bbZ$.

For each $v \in \bbC$ let $G_v$ be the isotropy subgroup $\{g \in G:\ \xi^gv = v\}$.  Then $G_v \cong \bbR$, and $G_v$ and $G_w$ are conjugated by the `translational' element $(w - v,0) \in G$.  Moreover, since any translation of $\bbC$ may be obtained as a composite of two rotations about different points, the groups $G_v$ together generate $G$, and so $G_v^\rm{n} = G$ for every $v$. A simple calculation shows that in coordinates one has
\[G_v = \{(v - \rho^t(v),t):\ t \in \bbR\}.\]

Now consider the action $\pi:G\actson L_\bbC^2(m_{\rm{S}^1})\cong L^2(m_{\rm{S}^1})\otimes_\bbR \bbC$ defined by
\[(\pi(u,t)f)(z) := \rm{e}^{2\pi\rm{i}\langle \rho^{-t}u,z\rangle}f(\rho^tz),\]
where $\langle \rho^{-t}u,z\rangle$ is the usual inner product of $\bbC$ regarded as a vector space over $\bbR$. (A routine check shows that this formula correctly defines an action of $G$.) The subspace $\rm{Fix}(\pi(G_v))$ consists of those functions $f$ such that
\[\rm{e}^{2\pi\rm{i}\langle \rho^{-t}u,z\rangle}f(\rho^tz) = f(z)\quad\forall z \in \rm{S}^1,\,t\in\bbR:\]
that is, of the constant complex multiples of the function $z\mapsto \rm{e}^{-2\pi\rm{i}\langle u,z\rangle}$. These are all distinct $2$-real-dimensional subspaces of $L_\bbC^2(m_{\rm{S}^1})$, so are not equal to $\rm{Fix}(\pi(G)) = \{0\}$, and also (by considering close-by values of $v$, for instance) they are not pairwise orthogonal. \fin

Another useful result in a similar vein to Corollary~\ref{cor:normal-clos} is the following simple relative of the Pugh-Shub Theorem~\cite{PugShu71}.  An adaptation of their theorem to the setting of nilpotent groups has previously been given by Ratner in Proposition 5.1 of~\cite{Rat91-a}.  Although our formulation is superficially different from hers, each version can easily be deduced from the proof of the other.

\begin{lem}\label{lem:PS}
Let $\pi:G\actson \frH$ be an orthogonal representation of a connected nilpotent Lie group, and let $\Lat\frg$ be the family of all proper Lie subalgebras of $\frg$.  Then the subfamily
\[\cal{A} := \{V \in \Lat\frg:\ \rm{Fix}(\pi(\exp V))\supsetneqq \rm{Fix}(\pi(G))\}\]
has countably many maximal elements.
\end{lem}

\textbf{Proof}\quad Suppose that $V_1,V_2\in \cal{A}$ are two distinct maximal elements.  Then the Lie subalgebra generated by $V_1 + V_2$ must strictly contain them both, and hence
\[\rm{Fix}(\pi(\langle \exp V_1\cup \exp V_2\rangle)) = \rm{Fix}(\pi(G)),\]
by their maximality.

Corollary~\ref{cor:rel-ind-over-common} now implies that $\rm{Fix}(\pi(\exp V_1))$ and $\rm{Fix}(\pi(\exp V_2))$ are relatively orthogonal over $\rm{Fix}(\pi(G))$.  Therefore there can be at most countably many of these maximal elements of $\cal{A}$, because $\frH$ is separable: indeed, if $\A_1\subseteq \A$ were an uncountable collection of maximal elements, then choosing some representative unit vectors
\[x_V \in \rm{Fix}(\pi(\exp V))\ominus \rm{Fix}(\pi(G))\quad\forall V \in \cal{A}_1\]
would give an uncountable sequence of orthonormal vectors in $\frH$, and hence a contradiction. \qed

\textbf{Example}\quad It is certainly not true that $\cal{A}_1$ is generally finite.  For example, consider the obvious rotation action of $\bbR^2$ on $\bbT^2$ and let $\pi:\bbR^2\actson L^2(m_{\bbT^2})$ be the resulting orthogonal representation.  Then any one-dimensional subgroup $\bbR\bf{v} \leq \bbR^2$ of rational slope has some non-trivial invariant functions, but the whole $\bbR^2$-action is ergodic. \fin

This conclusion of countability (rather than finitude) gives rise to the need for the notion of Zariski genericity (rather than simply Zariski openness).  The connection between them is established by the following.

\begin{cor}\label{cor:generically-const-fpspace}
If $\phi:\bbR\times \bbR^r\to G$ is a polynomial map into a connected and simply connected nilpotent Lie group and $\pi:G\actson \frH$ is an orthogonal representation, then the map
\[\bbR^r\to (\hbox{subspaces of $\frH$}):h\mapsto \rm{Fix}(\pi(\langle \img\phi(\cdot,h) \rangle))\]
takes the fixed value $\rm{Fix}(\pi(\langle\img\phi\rangle))$ Zariski generically.  Similarly, if $(X,\S,\mu,u)$ is a $G$-system then the $\s$-subalgebra $\S^{\langle \img \phi(\cdot,h)\rangle}$ agrees with $\S^{\langle \img \phi \rangle}$ up to $\mu$-negligible sets for Zariski generic $h$.
\end{cor}

\textbf{Proof}\quad Replacing $G$ with $\langle \img \phi\rangle^{\rm{n}}$ if necessary, we may assume they are equal.

Let $\cal{A} \leq \Lat \frg$ be the family of all Lie subalgebras with fixed-point subspaces strictly larger than $\rm{Fix}(\pi(G))$, as in Lemma~\ref{lem:PS}, and let $\cal{A}_1 \subseteq \cal{A}$ be the subfamily of maximal elements of $\cal{A}$, so Lemma~\ref{lem:PS} shows that this is countable.  Since $\rm{Fix}(\pi(\exp V^\rm{n})) = \rm{Fix}(\pi(\exp V))$ for any $V \in \Lat\frg$ by Corollary~\ref{cor:normal-clos}, by maximality we must have $V = V^\rm{n}$ for every $V \in \cal{A}_1$.

Now,
\begin{multline*}
\{h:\ \rm{Fix}(\pi(\langle \img \phi(\cdot,h) \rangle)) \supsetneqq \rm{Fix}(\pi(\langle \img \phi\rangle)) = \rm{Fix}(\pi(G))\}\\ = \bigcup_{V \in \cal{A}_1}\{h:\ \phi(t,h) \in \exp V\ \forall t\in\bbR\},
\end{multline*}
and so by the countability of $\cal{A}_1$ it suffices to show that each individual set $\{h:\ \exp\phi(t,h) \in V\ \forall t\in\bbR\}$ is proper and Zariski closed in $\bbR^r$. Since $\rm{Fix}(\pi(\langle \img \phi\rangle)) = \rm{Fix}(\pi(G))$, the subgroup $\langle \img \phi\rangle$ is not contained in $\exp V$ for any $V \in \cal{A}_1$, and so in fact $\img \phi \not\subseteq \exp V$ (since $\exp V$ is itself a subgroup).

Therefore for any $V \in \cal{A}_1$ we may choose a linear form $\ell \in \frg^\ast$ which annihilates $V$ but does not annihilate the whole of $\exp^{-1}\langle\img \phi\rangle$, and now one has
\[\{h:\ \phi(t,h) \in \exp V\ \forall t \in \bbR\} \subseteq \{h:\ \ell(\exp^{-1}(\phi(t,h))) = 0\ \forall t \in \bbR\}.\]

However, the map $(t,h)\mapsto \ell(\exp^{-1}(\phi(t,h)))$ is a polynomial $\bbR\times \bbR^n\to \bbR$, by Corollary~\ref{cor:exppoly-is-poly}.  By collecting monomials it may be expressed as
\[t^d p_d(h) + t^{d-1}p_{d-1}(h) + \cdots + t p_1(h) + p_0(h)\]
for some $p_i \in \bbR[h_1,\ldots,h_r]$, and now
\[\{h:\ \ell(\exp^{-1}(\phi(t,h))) = 0\ \forall t \in \bbR\} = \bigcap_{i=0}^d \{h:\ p_i(h) = 0\}.\]
This is manifestly a real algebraic subvariety of $\bbR^n$, and it is proper because the map $\ell\circ \exp^{-1}\circ \phi$ was chosen so as not to vanish identically, so it is a Zariski meagre subset of $\bbR^n$, as required.

Once again the conclusion about $G$-systems follows at once by considering Koopman representations. \qed

\section{Idempotent classes}\label{sec:idem}

The final ingredients needed for the proof of Theorem~\ref{thm:bigmain} are some results on `idempotent classes' of probability-preserving systems.  These were introduced in~\cite{Aus--lindeppleasant1,Aus--lindeppleasant2} building on the earlier notion of a `pleasant extensions' of systems~\cite{Aus--nonconv} (and also worth comparing with Host's `magic extensions' from~\cite{Hos09}).

\begin{dfn}[Idempotent and hereditary classes]
For any l.c.s.c. group $G$, a class $\sf{C}$ of jointly-measurable, probability-preserving $G$-systems is \textbf{idempotent} if it is closed under measure-theoretic isomorphisms, inverse limits and arbitrary joinings.  It is \textbf{hereditary} if it is closed under passing to factors.
\end{dfn}

\textbf{Example}\quad 
The leading examples of idempotent classes are those of the form
\[\sfC_0^{H_1}\vee \cdots \vee \sfC_0^{H_\ell}\]
for some closed normal subgroups $H_1,H_2,\ldots,H_\ell \unlhd G$, where this denotes the class of all $G$-systems which can be expressed as a joining of systems $\bfY_1$, $\bfY_2$, \ldots, $\bfY_\ell$ where each $\bfY_i$ has trivial $H_i$-subaction. \fin

The reference~\cite{Aus--thesis} contains an introduction to idempotent classes in the case of a discrete acting group.  In earlier works, idempotent classes were introduced to set up the theory of `sated extensions' of probability-preserving systems, which then play the primary r\^ole in applications of these ideas.  However, sated extensions are a little inconvenient in the present setting, and so we will work instead with some more elementary results about idempotent classes.  The reasoning behind this change of perspective relates to the need to change the group that acts on a system, which will appear in Section~\ref{sec:char-factor}.

In addition, our interest here is in actions of Lie groups, for which these ideas have not previously appeared in the literature.  Therefore the basic definitions and results we need have been included below for completeness.  Only very simple changes and additions are needed to the treatments in~\cite{Aus--thesis} or~\cite{Aus--lindeppleasant1}.  We will also introduce a slightly novel example of an idempotent class, useful for handling the polynomial maps of the present setting.

\begin{lem}[C.f. Lemma 2.2.2 in~\cite{Aus--thesis}]
If $\sfC$ is an idempotent class of $G$-systems and $\bfX = (X,\S,\mu,u)$ is any $G$-system, then $\bfX$ has an essentially unique largest factor $\L \leq \S$ that may be generated by a factor map to a member of $\sfC$.
\end{lem}

\textbf{Proof}\quad It is clear that under the above assumption the family of factors
\[\{\Xi\leq \S:\ \Xi\ \hbox{is generated by a factor map to a system in}\  \sfC\}\]	
is nonempty (it contains $\{\emptyset,X\}$, which corresponds to the trivial system), upwards directed (because $\sfC$ is closed under joinings) and closed under taking $\s$-algebra completions of increasing unions (because $\sfC$ is closed under inverse
limits). There is therefore a maximal $\s$-subalgebra in this family. \qed

\begin{dfn}[Maximal $\sfC$-factors]
The factor $\L$ obtained in the preceding lemma is the \textbf{maximal $\sfC$-factor} of $(X,\S,\mu,u)$, and will sometimes be denoted by the (slightly abusive) notation $\sfC\S$.  Similarly, we will sometimes denote by $\sfC\bfX$ a choice of a member of $\sfC$ such that $\sfC\S$ can be generated by a factor map $\bfX\to \sfC\bfX$.
\end{dfn}

The importance of idempotent classes derives from the following proposition.

\begin{prop}[Joinings to members of idempotent classes]\label{prop:idem-inverse}
Suppose that $\sfC$ is a hereditary idempotent class of $G$-systems, that $\bfX = (X,\S,\mu,u)$ is any $G$-system and $\bfY = (Y,\Phi,\nu,v)$ is a member of $\sfC$.  Then for any joining
\begin{center}
$\phantom{i}$\xymatrix{& \bfZ = (X\times Y,\S\otimes \Phi,\l,u\times v)\ar_{\pi}[dl]\ar^{\xi}[dr]\\
\bfX & &\bfY,
}
\end{center}
where $\pi$ and $\xi$ are the coordinate projections, there is some further factor $\L$ of $\bfX$ which is generated by a factor map to a member of $\sfC$ and such that the factor $\pi^{-1}(\S)$ is relatively independent from $\xi^{-1}(\Phi)$ over $\pi^{-1}(\L)$.  Concretely, this means that
\[\int_Z f(x)g(y)\,\l(\d x,\d y) = \int_Z \sfE_\mu(f\,|\,\L)(x)g(y)\,\l(\d x,\d y)\]
for any $f \in L^\infty(\mu)$ and $g \in L^\infty(\nu)$ (so we do not require that $\pi^{-1}(\L)$ also be contained in $\xi^{-1}(\Phi)$ up to negligible sets).
\end{prop}

\textbf{Proof}\quad We will construct from the joining $\l$ a new joining of $\bfX$ with a $\sfC$-system such that $\l$ is relatively independent over a factor of $\bfX$ which in that new joining is actually determined by the coordinate in the $\sfC$-system.

Let $\L:X\to \rm{Pr}\,Y$ be a disintegration of $\l$ over the coordinate projection to $X$. Form the infinite Cartesian product
\[Z' := X\times Y^\bbN\]
and let $\l'$ be the $(u\times v^{\times \bbN})$-invariant measure obtained as the relatively independent product of copies of $\l$:
\[\l' = \int_X \delta_x\otimes\L(x,\cdot)^{\otimes\bbN}\,\mu(\d x).\]
Let $\pi':Z'\to X$ be the first coordinate projection, and let $\l_1$ be the image of $\l'$ under the projection to $Y^\bbN$.

Finally, let $\L \leq \S$ be the $\s$-algebra of those sets which are $\l'$-a.s. determined by the remaining coordinates of $Z'$:
\[\L := \{A \in \S:\ \exists B \in \Phi^{\otimes \bbN}\ \rm{s.t.}\ \l'((A\times Y^\bbN)\triangle (X\times B)) = 0\}.\]
This is clearly a factor of $\bfX$, and by definition it also specifies a factor of the system $(Y^\bbN,\Phi^{\otimes\bbN},\l_1,v^{\times\bbN})$ (since each $A \in \L$ is identified with a member of $\Phi^{\otimes \bbN}$, uniquely up to negligible sets).  Let $\L' := (\pi')^{-1}(\L)$, so up to negligible sets this is measurable with respect to either $\pi'$ or the coordinate projection $Z'\to Y^\bbN$. The system $(Y^\bbN,\Phi^{\otimes\bbN},\l_1,v^{\times\bbN})$ is a member of $\sfC$, because $\bfY \in \sfC$ and $\sfC$ is closed under joinings; and hence the factor of $\bfX$ generated by $\L$ is also in $\sfC$, because it may be identified with a factor of that member of $\sfC$ and $\sfC$ is hereditary.

Now let $f \in L^\infty(\mu)$ and $g \in L^\infty(\nu)$.  To prove the desired equality of integrals, it suffices to show that
\[\sfE_\mu(f\,|\,\L) = 0 \quad \Longrightarrow\quad \sfE_\l(f\circ \pi\,|\,\{\emptyset,X\}\otimes \Phi) = 0,\]
since an arbitrary $f$ may be decomposed as $\sfE_\mu(f\,|\,\L) + (f - \sfE_\mu(f\,|\,\L))$, and this decomposition inserted into the two integrals against $g$ then shows that they are equal.

Thus, suppose conversely that
\[g := \sfE_\l(f\circ \pi\,|\,\{\emptyset,X\}\otimes \Phi) \neq 0,\]
and hence
\[\int_Z (f\circ \pi)\cdot g\,\d\l = \|g\|_2^2 \neq 0.\]

For each $i \in \bbN$ let $\a_i:Z'\to Y$ be the coordinate projection to the $i^\rm{th}$ copy of $Y$, and let $g_i := g\circ\a_i$.  By the construction of $\l'$, the pair of coordinates $(\pi,\a_i):Z'\to Z$ has the distribution $\l$ for any $i$. This has the following two consequences:
\begin{itemize}
\item for any $M\geq 1$ one has
\[\int_{Z'}(f\circ \pi')\Big(\frac{1}{M}\sum_{m=1}^Mg_m\Big)\,\d\l' = \int_Z (f\circ \pi)\cdot g\,\d\l = \|g\|_2^2 \neq 0;\]
\item for all $i$ one has
\[\sfE_{\l'}(g_i\,|\,\S\otimes \{\emptyset,Y^\bbN\}) = \sfE_{\l'}(g_1\,|\,\S\otimes \{\emptyset,Y^\bbN\}),\]
so we may let $h$ be this common conditional expectation.
\end{itemize}

Next, since all the $Y$-valued coordinates in $Z'$ are relatively independent under $\l'$ given the $X$-coordinate, one has
\[\int_{Z'}(g_i - h)(g_j - h)\,\d\l' = 0 \quad\hbox{whenever}\ i\neq j,\]
and as $M\to\infty$ this implies the simple estimate
\[\Big\|\frac{1}{M}\sum_{m=1}^Mg_m - h\Big\|_2^2 = \Big\|\frac{1}{M}\sum_{m=1}^M(g_m - h)\Big\|_2^2 = \frac{1}{M^2}\sum_{m=1}^M\|g_m - h\|_2^2 = \rm{O}\Big(\frac{1}{M}\Big).\]
Hence
\[\frac{1}{M}\sum_{m=1}^Mg_m \to h\]
in $\|\cdot\|_2$ as $M\to\infty$.  On the one hand, this implies that $h$ is a limit of functions measurable with respect to $\{\emptyset,X\}\otimes \Phi^{\otimes\bbN}$, hence is itself virtually measurable with respect to that $\s$-algebra.  Therefore as a function on $X$ it must actually be $\L$-measurable.  On the other hand, the above non-vanishing integral now gives
\[\int_{Z'}(f\circ \pi')\cdot h\,\d\l' \neq 0.\]
Therefore $\sfE_\mu(f\,|\,\L)\neq 0$, so since $\L$ defines a $\sfC$-factor of $\bfX$ this completes the proof. \qed

\textbf{Remark}\quad This proof can be presented in several superficially different ways.  On the one hand, it can be deduced almost immediately from a well-chosen appeal to the de Finetti-Hewitt-Savage Theorem, as in the paper~\cite{LesRitdelaRue03} of Lesigne, Rittaud and de la Rue (see also Section 8.5 in Glasner~\cite{Gla03}).  On the other, it is a close cousin of the proof that for any idempotent class $\sfC$, any system $\bfX$ has an extension that is `$\sfC$-sated' (Theorem 2.3.2 in~\cite{Aus--thesis}).  \fin

%\textbf{Remark}\quad One can think of satedness as a distant relative of injectivity for a module $M$ over a ring $R$.  Injectivity asserts that given an inclusion of $R$-modules $A \into B$ and a homomorphism $h:A\to M$, the homomorphism can be lifted to $B$: equivalently, it can be factorized as $A \into B\to M$.  Insofar as the obstructions to finding such an extended homomorphism are the linear equations over $R$ that hold among a set of generators of $B$ and arbitrary elements of $A$, the injectivity of $M$ asserts that for any such system of $R$-linear equations a solution can be found in $M$ with the relevant elements of $A$ replaced by their images in $h(A)$.  Satedness may be phrased similarly, except that it applies to finding statistical correlations rather than solutions to equations.  If $\bfX$ is $\sfC$-sated, and we know that for some $f \in L^\infty(\mu)$ and extension $\pi: \t{\bfX}\to\bfX$ there is a non-trivial conditional expectation $\sfE(f\circ\pi\,|\,\sfC\t{\S}) \neq 0$, then we can find a $\sfC$-factor of the original system $\bfX$ which also witnesses this. A slight extension of this last principle will be given formally as Lemma~\ref{lem:corn-descent} below, in preparation for future use. \fin

In previous applications, the idempotent classes of importance were those of the form $\sfC_0^{H_1}\vee\cdots\vee \sfC_0^{H_\ell}$, introduced as examples above.  Here we will need some slightly more complicated examples, because in order to account for the possible relations among the polynomials of a tuple $\F$ we will need to consider simultaneously actions of $G$ and also some `more free' covering group $q:\t{G}\to G$.

\begin{lem}\label{lem:still-idemp-1}
Suppose that $q:H\to G$ is a continuous homomorphism of l.c.s.c. groups and that $\sfC$ is an idempotent class of $H$-systems.  Then
\[q_\ast\sfC := \{\hbox{$G$-systems $\bfX$ such that $\bfX^{q(\cdot)} \in \sfC$}\}\]
is an idempotent class of $G$-systems, and it is hereditary if $\sfC$ is hereditary.
\end{lem}

\textbf{Proof}\quad We must verify that $q_\ast\sfC$ is closed under joinings and inverse limits.  Both are immediate: if $\bfY$ is a joining of $\bfX_i \in q_\ast\sfC$ for $i=1,2$ then $\bfY^{q(\cdot)}$ is the corresponding joining of $\bfX_i^{q(\cdot)}$, so lies in $\sfC$ because $\sfC$ is closed under joinings, and similarly for inverse limits. The last assertion also follows at once from the definition. \qed

\begin{dfn}
The new class $q_\ast\sfC$ constructed in the previous lemma is the \textbf{image} of $\sfC$ under $q$.
\end{dfn}

\begin{lem}\label{lem:still-idemp-2}
If $\sfC$ is an idempotent class of $G$-systems then
\[\hat{\sfC} := \{\bfX:\ \bfX\ \hbox{is a factor of a member of}\ \sfC\}\]
is a hereditary idempotent class.
\end{lem}

\textbf{Proof}\quad The hereditary property is built into the definition, so once again it remains to check closure under joinings and inverse limits.  Both are routine, so we give the proof only for joinings.  Suppose that $\bfY_i = (Y_i,\Phi_i,\nu_i,v_i) \in \hat{\sfC}$ for $i=1,2$, that $\pi_i:\bfX_i\to \bfY_i$ are factors with $\bfX_i = (X_i,\S_i,\mu_i,u_i)\in \sfC$ for $i=1,2$, and that $\bfZ = (Y_1\times Y_2, \Phi_1 \otimes \Phi_2, \l, v_1\times v_2)$ defines a joining of $\bfY_1$ and $\bfY_2$. Then we may define a joining of $\bfX_1$ and $\bfX_2$ as a relatively independent product: letting $P_i:Y_i\to \Pr(X_i)$ be a probability kernel representing the disintegration of $\mu_i$ over $\pi_i$, define
\[\l':= \int_{Y_1\times Y_2}P(y_1,\cdot)\otimes P(y_2,\cdot)\,\l(\d y_1,\d y_2).\]

Now $(X_1\times X_2,\S_1\otimes \S_2,\l',u_1\times u_2)$ is a joining of $\bfX_1$ and $\bfX_2$, and hence a member of $\sfC$. The map $(x_1,x_2)\mapsto (\pi_1(x_1),\pi_2(x_2))$ witnesses $\bfZ$ as a factor of this member of $\sfC$, so $\bfZ \in \hat{\sfC}$. \qed

\begin{dfn}\label{dfn:down-clos}
The class $\hat{\sfC}$ constructed above is the \textbf{downward closure} of $\sfC$.
\end{dfn}

When we come to apply this machinery, satedness relative only to classes of the form $\sfC_0^{H_1} \vee \cdots \vee \sfC_0^{H_\ell}$ will not give us quite enough purchase over our situation.  Instead we will need to first form an extended group $q:\t{G}\onto G$ (in which copies of certain subgroups of $G$ have been made `more independent': see Section~\ref{sec:char-factor}), and then for some subgroups $\t{H}_1$, $\t{H}_2$, \ldots, $\t{H}_\ell \unlhd \t{G}$ we will need to use satedness relative to the class
\[q_\ast\big(\ (\sfC_0^{\t{H}_1} \vee \cdots \vee \sfC_0^{\t{H}_\ell})^\wedge\ \big).\]
In prose, this is
\begin{quote}
`The class of $G$-systems which, upon re-writing them as $\t{G}$-systems, become factors of joinings of systems in which one of the $\t{H}_i$ acts trivially.'
\end{quote}
This manoeuvre will appear during the proof of Proposition~\ref{prop:vdC-appn} below, where the need for it will become clearer.  The particular way in which we will appeal to satedness with respect to such a class is captured by the following lemma.

\begin{lem}\label{lem:corn-descent}
Suppose that $q:H\onto G$ is a continuous epimorphism of Lie groups, that $\sfC$ is an idempotent class of $H$-systems and that $\bfX = (X,\S,\mu,u)$ is a $G$-system.  In addition, suppose that $f \in L^\infty(\mu)$ and that
\[\pi:\bfY = (Y,\Phi,\nu,v)\to \bfX^{q(\cdot)}\] is an extension of $H$-systems such that \[\sfE_\nu(f\circ \pi\,|\,\sfC\Phi) \neq 0.\]
Then also
\[\sfE_\mu(f\,|\,(q_\ast\hat{\sfC})\S) \neq 0.\]
\end{lem}

\textbf{Proof}\quad We have $\sfE_\nu(f\circ \pi\,|\,\sfC\Phi) \neq 0$ by assumption, but on the other hand the function $f\circ\pi$ is invariant under $v^h$ for every $h \in \ker q$:
\[f\circ\pi\circ v^h = f\circ u^{q(h)}\circ \pi = f\circ u^e\circ\pi = f\circ \pi.\]
Since $\sfC\Phi$ is a factor of the whole $H$-action $v$, the conditional expectation operator $\sfE_\nu(\,\cdot\,|\,\sfC\Phi)$ preserves this $\ker q$-invariance.  Therefore $\sfE_\nu(f\circ \pi\,|\,\sfC\Phi)$ is measurable not only with respect to $\sfC\Phi$ but also with respect to $\Phi^{\ker q}$.

Let $\a:\bfY\to \bfZ$ be a factor map onto another system which generates the factor $\Phi^{\ker q}\cap \sfC\Phi \leq \Phi$, so its target system $\bfZ$ is an element of $\hat{\sfC}$ and has $\ker q$ acting trivially.  Therefore this action of $H$ may be identified with an action of $G$ composed through $q$, say $\bfZ = \bfW^{q(\cdot)}$ for some $G$-system $\bfW$.  (The joint measurability of $v$ implies that of the action of $G$ on $\bfW$, simply by choosing an everywhere-defined Borel selector $G\to H$, as we clearly may for Lie group epimorphisms because they are are locally diffeomorphic to orthogonal projections.)

Now the diagram
\begin{center}
$\phantom{i}$\xymatrix{
& \bfY\ar[dl]_\pi\ar[dr]^\a\\
\bfX^{q(\cdot)} & & \bfW^{q(\cdot)}
}
\end{center}
defines a joining of $\bfX^{q(\cdot)}$ and $\bfW^{q(\cdot)}$.  It therefore also defines a joining of $\bfX$ and $\bfW$, by simply identifying it with an invariant measure on $X\times W$ and writing the actions in terms of $G$ rather than $H$.

Our assumption on $f$ gives that $\sfE(f\circ \pi\,|\,\a) \neq 0$.  Therefore, within this joining of $\bfX$ and $\bfW$, the lift of $f$ has non-trivial conditional expectation onto the copy of $\bfW$, which is a member of $q_\ast\hat{\sfC}$, and so by Proposition~\ref{prop:idem-inverse} and Lemma~\ref{lem:still-idemp-2} this implies $\sfE_\mu(f\,|\,(q_\ast\hat{\sfC})\S) \neq 0$. \qed

\section{The case of two-fold joinings}\label{sec:k=2}

The case of Theorem~\ref{thm:main} in which $k=1$ will form the base of an inductive proof of the full theorem, and must be handled separately.  Its proof is quite routine in the shadow of other works in this area, but it does already contain an appeal to the van der Corput estimate and an induction on the PET ordering for single polynomials (rather than whole tuples).  It thus serves as a helpful preparation for the full induction that is to come.

\begin{prop}\label{prop:k=1}
Suppose that $\pi:G\actson \frH$ is an orthogonal representation and $\phi:\bbR\times \bbR^r \to G$ is a polynomial map such that $\phi(0,\cdot) \equiv e$.  Then the operator averages
\[\barint_0^T \pi(\phi(t,h))\,\d t\]
converge in the strong operator topology for every $h$, and the limit operator $P_h$ is Zariski generically equal to the orthoprojection onto $\rm{Fix}(\pi(\langle \img\phi\rangle))$.
\end{prop}

\textbf{Proof}\quad\textbf{Step 1}\quad First suppose that $\phi$ is linear in the first coordinate, meaning that $\phi(\cdot,h)$ is a homomorphism for every $h \in \bbR^r$.  Then for every $h$ the map $t\mapsto \phi(t,h)$ takes values in a $1$-parameter subgroup of $G$, and so the classical ergodic theorem for orthogonal flows gives
\[\barint_0^T \pi(\phi(t,h))\,\d t\stackrel{\rm{SOT}}{\to} P_h,\]
where $P_h$ is the orthoprojection onto $\rm{Fix}(\pi(\langle \img\phi(\cdot,h)\rangle))$.  By Corollary~\ref{cor:generically-const-fpspace} this equals $\rm{Fix}(\pi(\langle \img \phi\rangle))$ Zariski generically, and so the proof is complete in the linear case.

\quad\textbf{Step 2}\quad For arbitrary polynomial maps $\phi$ we show by PET induction that if
\[\barint_0^T \pi(\phi(t,h))v\,\d t\,\,\not\!\!\to 0\]
for some $v \in \frH$, then $P_hv \neq 0$, where again $P_h$ is the orthoprojection onto $\rm{Fix}(\pi(\langle \img \phi(\cdot,h)\rangle))$.  By decomposing an arbitrary $v$ as $(1 - P_h)v + P_hv$ and appealing to Corollary~\ref{cor:generically-const-fpspace} again, this will complete the proof.

If
\[\barint_0^T \pi(\phi(t,h))v\,\d t\,\,\not\!\!\to 0\]
then the van der Corput estimate~\ref{lem:vdC} gives that also
\begin{multline*}
\barint_0^S\barint_0^T \langle\pi(\phi(t+s,h))v,\pi(\phi(t,h)v\rangle\,\d t\,\d s\\
=\Big\langle\barint_0^S\barint_0^T \pi(\phi(t,h)^{-1}\phi(t+s,h))v\,\d t\,\d s,\ v\Big\rangle \,\,\not\!\!\to 0
\end{multline*}
as $T\to\infty$ and then $S\to \infty$.

By the special case of Lemma~\ref{lem:PET-calcns} for singleton families we have
\[\{(t,s,h)\mapsto \phi(t,h)^{-1}\phi(t+s,h)\}\prec_\PET\{\phi\},\]
and so the inductive hypothesis gives
\[\barint_0^T \pi(\phi(t,h)^{-1}\phi(t+s,h))v\,\d t \to Q_{s,h}v \quad\quad \hbox{as}\ T\to\infty\]
with $Q_{s,h}$ the orthoprojection onto $\rm{Fix}(\pi(\langle \img\phi(\cdot,h)^{-1}\phi(\cdot+s,h)\rangle ))$.

By Corollary~\ref{cor:generically-const-fpspace}, for every fixed $h$ we have
\[\rm{Fix}(\pi(\langle \img\phi(\cdot,h)^{-1}\phi(\cdot+s,h)\rangle )) = \rm{Fix}(\pi(\langle \img\phi(\cdot,h)^{-1}\phi(\cdot+ \cdot,h)\rangle ))\]
for Zariski generic $s$, and now since $\phi(0,h) \equiv e$ this is equal to
\[\rm{Fix}(\pi(\langle\img \phi(\cdot,h)\rangle)).\]
In particular, for every $h$ this equality must hold for Lebesgue-a.e. $s$, and thus our previous average over $s$ may be written instead as
\[\barint_0^S Q_{s,h}v\,\d s = \barint_0^S P_hv\, \d s \equiv P_h v.\]
This proves that $P_hv \neq 0$, as required. \qed

\section{A partially characteristic factor}\label{sec:char-factor}

Now fix the following assumptions for this section and the next:
\begin{itemize}
\item $G$ is an $s$-step connected and simply connected nilpotent Lie group;
\item $\F = (\phi_1,\phi_2,\ldots,\phi_k)$ is a tuple of polynomial maps $\bbR\times \bbR^r\to G$ with $k\geq 2$ in which $\phi_1$ is a pivot, such that $\phi_i(0,\cdot) \equiv e$ for each $i$, and such that $G = \langle\img \phi_1\cup\cdots\cup \img\phi_k\rangle$ (otherwise we may simply replace $G$ with this smaller group);
\item $(X_i,\S_i,\mu_i,u_i)$ for $0 \leq i \leq k$ is a tuple of $G$-systems, and $\l$ is a joining of them;
\item $A^\l_T$ for $T \in [0,\infty)$ is the family of averaging operators associated to the orbit of $\l$ under $(\phi_1,\phi_2,\ldots,\phi_k)$ as in Theorem~\ref{thm:bigmain}, so note that these implicitly depend on $h$, the parameter in the argument of the $\phi_i$ which is \emph{not} averaged.
\end{itemize}

At the heart of the inductive proof of Theorem~\ref{thm:bigmain} lies a result promising that in order to study the functional averages $A^\l_T(f_1,f_2,\ldots,f_k)$, one may assume that one of the functions $f_i$ has some special additional structure (which we will see later enables a further reduction to the case of a simpler family of polynomial maps).  This extra structure is captured by a simple adaptation of an important idea introduced in~\cite{FurWei96}, and which has been used extensively since (see, for instance,~\cite{HosKra05,Zie07,Aus--nonconv,Aus--lindeppleasant1}).

\begin{dfn}[Partially characteristic factor]
In the above setting a factor $\L \leq \S_1$ is \textbf{partially characteristic} for the averages $A^\l_T$ if for any tuple of functions $f_i \in L^\infty(\mu_i)$ one has
\[\big\|A^\l_T(f_1,f_2,\ldots,f_k) - A^\l_T\big(\sfE(f_1\,|\,\L),f_2,\ldots,f_k\big)\big\|_2 \to 0\]
as $T\to\infty$ for Zariski generic $h$ (recalling that the operators $A^\l_T$ implicitly depend on $h \in \bbR^r$).
\end{dfn}

\textbf{Remark}\quad The main difference between this definition and its predecessors in earlier papers is that here, in consonance with the statement of Theorem~\ref{thm:bigmain}, we require convergence only for Zariski generic $h$.

As stated, this definition allows the Zariski meagre set $F\subseteq \bbR^r$ containing those $h$ for which convergence fails to depend on $f_1$, $f_2$, \ldots, $f_k$.  However, it is easily checked that for a given $h$, this convergence holds for all tuples of functions if one knows that it holds for tuples drawn from some $\|\cdot\|_2$-dense subsets of the unit balls of $L^\infty(\mu_i)$, $i=1,2,\ldots,k$.  Since one can choose countable such subsets, we deduce that there is a countable intersection of Zariski residual subsets of $\bbR^r$ (which is therefore still Zariski residual) on which the above convergence holds for all tuples of functions.  \fin

As in many of the earlier works cited above, the first step towards proving the convergence of $A^\l_T(f_1,\ldots,f_k)$ will be to identify a partially characteristic factor with some useful structure. However, a new twist appears in the present setting: here we must first pass from $G$-systems to actions of some covering group of $G$.

To be precise, let
\[\t{\phi}_1:(t,h)\mapsto (\phi_1(t,h),\ldots,\phi_k(t,h)),\]
let
\[\t{\phi}_i:(t,h)\mapsto (\phi_i(t,h),\ldots,\phi_i(t,h)),\quad \hbox{for}\ i=2,3,\ldots,k,\]
(notice the subscripts in different coordinates), and let
\[\t{G} := \langle \img\t{\phi}_1\ \cup\ \img\t{\phi}_2\ \cup\ \cdots\ \cup\ \img\t{\phi}_k\rangle \leq G^{k+1}.\]
Let $q:\t{G} \to G$ be the restriction to $\t{G}$ of the projection $G^k\to G$ onto the first coordinate.  Then $q$ intertwines each $\t{\phi}_i$ with $\phi_i$ for $i\geq 1$ (because $\phi_i$ appears in the first coordinate of $\t{\phi}_i$ for every $i$).

It is easy to verify that $q(\t{G}) = G$.  The group $\t{G}$ is connected, because each $\t{\phi}_i(\cdot,h)$ passes through the origin for every $h$, and hence $\t{G} = \exp V$ for some Lie subalgebra $V \leq \frg^k$.  The image of $V$ under the first coordinate projection is a Lie subalgebra $V_1 \leq \frg$, and since $G$ is simply connected it follows that $\exp V_1$ is a closed subgroup of $G$ which is contained in $q(\t{G})$.  On the other hand it must contain $\img \phi_i$ for every $i\leq k$, so in fact $q(\t{G}) = \exp V_1 = G$.

The next technical proposition lies at the heart of all that follows. It provides a partially characteristic factor of $\bfX_1 = (X_1,\S_1,\mu_1,u_1)$ for the averages $A^\l_T$, but only at the cost of regarding instead the modified system $\bfX_1^{q(\cdot)}$.  The need for this sleight of hand will become clear during the proof.

\begin{prop}\label{prop:vdC-appn}
Assume that conclusions (1--3) of Theorem~\ref{thm:bigmain} have already been established for all polynomial families preceding $\F$ in the PET ordering, suppose that $\phi_1$ is a pivot, and let
\[\sfC := q_\ast\Big(\ \Big(\sfC_0^{\langle\img\t{\phi}_1\rangle}\vee\bigvee_{j=2}^k\sfC_0^{\langle\img\t{\phi}_j\t{\phi}_1^{-1}\rangle}\Big)^\wedge\ \Big).\]
(Recall the discussion following Definition~\ref{dfn:down-clos}.) Then for any systems $\bfX_i$, $i=0,1,\ldots,k$, the factor $\sfC\S_1 \leq \S_1$ is partially characteristic.
\end{prop}

\textbf{Remark}\quad Of course, once this proposition has been proved then it implies some conclusion even if $A^\l_T(f_1,\ldots,f_k) \,\,\not\!\!\to 0$ for just one value of $h$, because by fixing that $h$ we may simply regard each $\phi_i$ as a polynomial function of $t$ alone, and so apply the proposition with $r=0$.  Indeed, we will use this trick a few times later. However, one must beware of the delicacy that the idempotent class appearing in this proposition may not be the same after one makes such a restriction, so nor will the $\s$-sigma algebra $\sfC\S_1$ in general.  Even the group extension $q:\t{G}\to G$ itself will not be the same as above, but will depend on the choice of $h$.  Since at some points later we will really need the above conclusion about the generic behaviour of the averages in $h$, it seems easiest to formulate it as here and then apply it with a restricted parameter space when convenient. \fin

\textbf{Proof}\quad Since any $f_1$ may be decomposed as
\[\sfE_\mu(f_1\,|\,\sfC\S) + \big(f_1 - \sfE_\mu(f_1\,|\,\sfC\S)\big)\]
and the operator $A^\l_T$ is multilinear, it is enough to prove that if $\sfE_\mu(f_1\,|\,\sfC\S) = 0$ then for any $f_2$, \ldots, $f_k$ one has
\[\|A^\l_T(f_1,f_2,\ldots,f_k)\|_2\to 0\]
as $T\to\infty$ for Zariski generic $h$.  Contrapositively, this is equivalent to showing that if the set
\[E := \{h\in \bbR^r:\ \|A^\l_T(f_1,f_2,\ldots,f_k)\|_2 \,\,\not\!\!\to 0\ \hbox{as}\ T\to\infty\}\]
is not Zariski meagre then $\sfE_\mu(f_1\,|\,\sfC\S)\neq 0$.  Henceforth we assume that $E$ is not Zariski meagre.

Furthermore, in view of Lemma~\ref{lem:corn-descent}, it now suffices to find an extension of spaces $\pi:(\t{X},\t{\S},\t{\mu})\to (X_1,\S_1,\mu_1)$ and an action $\t{u}:\t{G}\actson (\t{X},\t{\S},\t{\mu})$ such that $\pi\circ \t{u} = u^{q(\cdot)}$ and
\[\sfE(f_1\circ \pi\,|\,\L) \neq 0,\]
where
\[\L := \t{\S}^{\langle \img\t{\phi}_1\rangle}\vee \bigvee_{i=2}^k \t{\S}^{\langle \img\t{\phi}_1\cdot \t{\phi}_i^{-1}\rangle}.\]
This is the point at which we have made use of the general properties of idempotent classes.  This implication will follow in two steps: applying the van der Corput estimate (Lemma~\ref{lem:vdC}), and interpreting what it tells us.

\quad\textbf{Step 1}\quad Letting 
\[g_{t,h} := M^\l\big(f_1\circ u_1^{\phi_1(t,h)},f_2\circ u_2^{\phi_2(t,h)},\ldots,f_k\circ u_k^{\phi_k(t,h)}\big),\]
the van der Corput estimate implies that for $h \in E$ one also has
\[\barint_0^S\barint_0^T\int_{X_0}g_{t+s,h}g_{t,h}\,\d\mu_0\,\d t\,\d s \,\,\not\!\!\to 0\]
as $T\to\infty$ and then $S\to\infty$.

For each $s$, by Lemma~\ref{lem:l_1l} we may re-write the two inner integrals here as
\begin{multline*}
\barint_0^T\int_{X_1^2\times \cdots \times X_k^2} (f_1\circ u_1^{\phi_1(t,h)})\otimes (f_1\circ u_1^{\phi_1(s,h)\psi_1(t,s,h)})  \otimes\\ \quad\quad\quad\quad\quad\quad\quad\quad \cdots \otimes (f_k \circ u_k^{\phi_k(t,h)})\otimes (f_k\circ u_k^{\phi_k(s,h)\psi_k(t,s,h)})\,\d(\l\otimes_0 \l)\,\d t,
\end{multline*}
where
\[\psi_i(t,s,h) := \phi_i(s,h)^{-1}\phi_i(t+s,h) \quad\quad \hbox{for each}\ i = 1,2,\ldots,k,\]
so $\psi_i:\bbR\times \bbR\times \bbR^r \to G$ is a polynomial map with the property that $\psi_i(0,\cdot,\cdot)\equiv e$.

Since $\l\otimes_0\l$ is a joining of two duplicates of each of the $G$-systems $(X_i,\S_i,\mu_i,u_i)$ for $1 \leq i \leq k$, it is invariant under the diagonal transformations $u_\Delta^{\varphi_1(t,h)^{-1}}$.  Applying this within the above integral shows that it is equal to
\begin{multline*}
\barint_0^T\int_{X_1^2\times \cdots \times X_k^2} f_1\otimes (f_1\circ u_1^{\phi_1(s,h)\psi'_1(t,s,h)}) \otimes\\ \quad\quad\quad\quad\quad\quad\quad\quad \cdots \otimes (f_k \circ u_k^{\phi'_k(t,h)})\otimes (f_k\circ u_k^{\phi_k(s,h)\psi'_k(t,s,h)})\,\d(\l\otimes_0 \l)\,\d t
\end{multline*}
with
\begin{eqnarray*}
\psi'_i(t,s,h) &:=& \psi_i(t,s,h)\phi_1(t,h)^{-1}\quad\hbox{for}\ i\geq 1\ \hbox{and}\\
\phi'_i(t,h) &:=& \phi_i(t,h)\phi_1(t,h)^{-1}\quad\hbox{for}\ i\geq 2.
\end{eqnarray*}
We recognize these as comprising the $1^{\rm{st}}$ derived family of $\F$, which by Lemma~\ref{lem:PET-calcns} precedes $\F$ in the PET ordering because $\phi_1$ was a pivot. Let
\[\vec{\psi}:(t,s,h) \mapsto (e,\psi'_1(t,s,h),\phi'_2(t,h),\psi'_2(t,s,h),\cdots,\phi_k'(t,h),\psi_k'(t,s,h)).\]

By the inductive hypothesis, for every $h \in \bbR^r$ there are a Zariski residual set $F_h \subseteq \bbR$ and a joining $\theta^h$ on $X_1^2 \times X_2^2\times\cdots \times X_k^2$ invariant under
\[\langle G^{\Delta 2k} \cup \img \vec{\psi}(\cdot,\cdot,h)\rangle\]
such that for all $s \in F_h$ the above integral tends to
\[\int_{X_1^2\times \cdots \times X_k^2} f_1\otimes (f_1\circ u_1^{\phi_1(s,h)})  \otimes \cdots \otimes f_k \otimes (f_k\circ u_k^{\phi_k(s,h)})\,\d\theta^h\]
as $T \to \infty$. Moreover these $\theta^h$ are equal to one fixed joining $\theta$ on a Zariski residual set of $h$, so that this $\theta$ must in fact be invariant under $\langle G^{\Delta 2k} \cup \img \vec{\psi} \rangle$.

Since the Zariski residual set $F_h$ has full Lebesgue measure, for each $h$ our previous average over $s$ may now be replaced by
\[\barint_0^S\int_{X_1^2\times \cdots \times X_k^2} f_1\otimes (f_1\circ u_1^{\phi_1(s,h)}) \otimes \cdots \otimes f_k \otimes (f_k\circ u_k^{\phi_k(s,h)})\,\d\theta^h\,\d s,\]
implying that for $h \in E$ this also does not vanish as $S \to \infty$.

Next, one has
\begin{eqnarray*}
&&(\phi_1(s,h),\phi_1(s,h),\ldots,\phi_k(s,h),\phi_k(s,h))\\
&&= (e,e,\ldots,\phi_k(s,h)\phi_1(s,h)^{-1},\phi_k(s,h)\phi_1(s,h)^{-1})\\
&&\quad\quad\quad\quad\quad\quad\quad \cdot (\phi_1(s,h),\phi_1(s,h),\ldots,\phi_1(s,h),\phi_1(s,h))\\
&&= \vec{\psi}(s,0,h)\cdot (\phi_1(s,h),\phi_1(s,h),\ldots,\phi_1(s,h),\phi_1(s,h))\\
&&\in \langle G^{\Delta 2k} \cup \img \vec{\psi}(\cdot,\cdot,h)\rangle
\end{eqnarray*}
for every $s$, and so each joining $\theta^h$ is already invariant under the new off-diagonal polynomial flow
\[\xi(\cdot,h):s \mapsto (\phi_1(s,h),\phi_1(s,h),\ldots,\phi_k(s,h),\phi_k(s,h)).\]

Since we may re-write the above average as
\[\barint_0^S\int_{X_1^2\times \cdots \times X_k^2} (f_1\otimes 1\otimes  \cdots \otimes f_k\otimes 1)\cdot \big((1\otimes f_1\otimes \cdots \otimes 1 \otimes f_k)\circ u_\times^{\xi(s,h)}\big)\,\d\theta^h\,\d s,\]
by the base case Proposition~\ref{prop:k=1} it must converge to
\[\int_{X_1^2\times \cdots \times X_k^2} (f_1\otimes 1\otimes  \cdots \otimes f_k\otimes 1)\cdot \sfE(1\otimes f_1\otimes \cdots \otimes 1 \otimes f_k\,|\,\S_\times^{\langle \img \xi(\cdot,h)\rangle})\,\d\theta^h\]
as $S\to\infty$, where $\S_\times := \S_1^{\otimes 2}\otimes \cdots \otimes \S_k^{\otimes 2}$ and the conditional expectation here is with respect to $\theta^h$.

Therefore this last integral is nonzero for every $h \in E$. Since the sets
\[\{h:\ \theta^h \neq \theta\}\quad\quad \hbox{and}\quad\quad \{h:\ \S_\times^{\langle \img \xi(\cdot,h)\rangle} \neq \S_\times^{\langle \img \xi\rangle}\ \hbox{up to $\theta$-negligible sets}\}\]
both \emph{are} Zariski meagre (the latter by Corollary~\ref{cor:generically-const-fpspace}), their union cannot contain $E$, and so any value $h \in E$ that is not in either of these meagre sets witnesses that
\[\int_{X_1^2\times \cdots \times X_k^2} (f_1\otimes 1\otimes  \cdots \otimes f_k\otimes 1)\cdot \sfE(1\otimes f_1\otimes \cdots \otimes 1 \otimes f_k\,|\,\S_\times^{\langle \img \xi \rangle})\,\d\theta \neq 0.\]

\quad\textbf{Step 2}\quad Now set
\[(\t{X},\t{\S},\t{\mu}) := \Big(\prod_{i=1}^kX_i^2,\bigotimes_{i=1}^k\S_i^{\otimes 2},\theta\Big)\]
and let $\pi:\t{X}\to X_1$ be the coordinate projection onto the first copy of $X_1$.  Observe that the polynomial map $\xi$ defined in Step 1 is simply a copy of $\t{\phi}_1$ in which each coordinate has been duplicated.  Define $q_1:\t{G}\to G^{2k}$ to be the restriction to $\t{G}$ of the coordinate-duplicating map
\[(g_1,g_2,\ldots,g_k)\mapsto (g_1,g_1,g_2,g_2,\ldots,g_k,g_k).\]
Composing $q_1$ with the Cartesian product action $u_\times$ of $G^{2k}$ now gives an action $\t{u}$ of $\t{G}$ on $(\t{X},\t{\S},\t{\mu})$, since we have already deduced from our inductive hypotheses that $\t{\mu} = \theta$ is invariant under $u_\Delta$ (and hence the image of $q_1\circ \t{\phi}_i$ for each $i\geq 2$) and also under $\langle\img \xi\rangle$ (which is the image of $q_1\circ\t{\phi}_1$).

On the first coordinate in $\prod_{i=1}^kX_i^2$, the transformation $\t{u}^g$ simply agrees with $u^g$ for any $g \in \langle\img\t{\phi}_2\cup\cdots\cup \img\t{\phi}_k\rangle$.  On the other hand,
\[\pi\circ \t{u}^{\t{\phi}_1(t,h)} \stackrel{\rm{def}}{=} \pi\circ (u_1^{\phi_1(t,h)}\times u_1^{\phi_1(t,h)}\times \cdots \times u_k^{\phi_k(t,h)}\times u_k^{\phi_k(t,h)}) = u_1^{\phi_1(t,h)}.\]
Since these cases together generate the whole of $\t{G}$, it follows that $\pi\circ \t{u}^{\t{g}} = u^{q(\t{g})}$ for all $\t{g} \in \t{G}$, where $q:\t{G}\to G$ is the covering homomorphism constructed previously.

Finally, an inspection of the action $\t{u}$ on the other coordinates of $\t{X}$ shows that
\begin{itemize}
\item for each $i\in \{2,3,\ldots,k\}$ the transformations $\t{u}^{\t{\phi}_1(t,h)}$ and $\t{u}^{\t{\phi}_i(t,h)}$ agree on the first coordinate copy of $X_i$, and
\item the function $\sfE(1\otimes f_1\otimes \cdots \otimes 1 \otimes f_k\,|\,\S_\times^{\langle \img \xi \rangle})$ is invariant under the $\t{u}$-action of $\langle\img\t{\phi}_1\rangle$.
\end{itemize}

Therefore the non-vanishing of the integral at the end of step 1 asserts that $f_1\circ\pi$ has a non-zero inner product with a function that is manifestly measurable with respect to a system in the class
\[\sfC_0^{\langle\img\t{\phi}_1\rangle}\vee\bigvee_{j=2}^k\sfC_0^{\langle\img\t{\phi}_j\t{\phi}_1^{-1}\rangle},\]
and hence $\sfE(f_1\circ\pi\,|\,\sfC\S_1)\neq 0$, as required. \qed

\textbf{Remarks}\quad\textbf{1.}\quad The above proof makes clear the need to extend the modified system $\bfX_1^{q(\cdot)}$, rather than $\bfX_1$ itself.  We constructed our extension from some joining on $X_1^2\times \cdots\times X_k^2$ through the coordinate projection onto $X_1$, and in order to derive the desired nonzero conditional expectation for it we needed the polynomial trajectory of transformations $u_1^{\phi_1(t,h)}$ downstairs to lift to the trajectory
\[u_1^{\phi_1(t,h)}\times u_1^{\phi_1(t,h)} \times u_2^{\phi_2(t,h)}\times u_2^{\phi_2(t,h)}\times \cdots \times u_k^{\phi_k(t,h)} \times u_k^{\phi_k(t,h)}.\]
The new map $\t{\phi}_1$ may not be a PET-minimal member of $(\t{\phi}_1,\ldots,\t{\phi}_k)$, and it also may not share its leading term with any of the lifts $\t{\phi}_i$ for $i\geq 2$, even if $\phi_1$ downstairs does have some leading terms in common with the other $\phi_i$.  Thus in order to write these $\t{\phi}_i$ as genuine lifts of the $\phi_i$ we must first split the group $G$ apart slightly in order to separate these leading terms.  Happily, the problem itself gives us a natural way to do this: the lifted polynomial mapping $\t{\phi}_1$ is suitably `separated' from $\t{\phi}_i$, $i\geq 2$, inside the Cartesian product $G^k$, so we have simply taken $\t{G}$ to be the closed subgroup of $G^k$ generated by these lifted mappings and composed our actions with the quotient map $q:\t{G}\to G$. 

\quad\textbf{2.}\quad If a factor $\L \leq \S_1$ is partially characteristic and we assume that the limits $\l^h = \lim_{T\to\infty}\l^h_T$ exist, then considering the integral formula
\[\int_{\prod_i X_i}f_0\otimes f_1\otimes \cdots \otimes f_k\,\d\l^h_T = \int_{X_0} f_0 \cdot A^\l_T(f_1,\ldots,f_k)\,\d\mu_0\]
shows that for Zariski generic $h$ the coordinate projection $\prod_i X_i\to X_1$ is relatively independent under $\l^h$ over its further factor generated by $\L \leq \S_1$.  Thus, knowledge of a non-trivial partially characteristic factor gives some structural information about the limit joining.

In particular, consider a case in which $q$ is an isomorphism (so that the subgroups $\langle\img \phi_1\rangle$ and $\langle \img\phi_2 \cup \cdots \cup \img \phi_k\rangle$ are already sufficiently `spread apart' in $G$), and suppose furthermore that the factor $\sfC\bfX_1$ can itself be expressed as a joining of systems $\bfZ_0 \in \sfC_0^{\langle \img \phi_1\rangle}$ and $\bfZ_i \in \sfC_0^{\langle \img(\phi_1\phi_i^{-1})\rangle}$ for $i\geq 2$ (rather than just as a factor of such).  Then we know that any limit joining $\l'$ must be relatively independent over the factor $\sfC\bfX_1$, and upon restricting ourselves to this factor we can express $\l'$ alternatively as a joining of
\[\bfX_0,\bfZ_0,\bfZ_2,\ldots,\bfZ_k,\bfX_2,\ldots,\bfX_k.\]
(In fact we will use a similar manipulation in the next section).  Moreover, the assumption that $\sfC\bfX_1$ itself be a joining is not terribly restrictive, since an arbitrary system $\bfX_1$ always has an extension for which this is true (by using the machinery of `$\sfC$-sated' extensions, as developed in Chapter 2 of~\cite{Aus--thesis}).

It would be interesting to know whether further use of the ideas behind Proposition~\ref{prop:vdC-appn} could give a more complete picture of the possible structure of $\l'$.  This would presumably involve repeated assertions of relative independence over increasingly `small' factors of the original system, on which increasingly large subgroups of $G$ act trivially. Such a picture does emerge in the study of the linear multiple averages constructed from a tuple of $\bbZ^d$-actions (see Chapter 4 of~\cite{Aus--thesis}), but in the present setting the need to keep track of a large family of different subgroups of $G$ may make the resulting description more obscure.

Even without a manageable description, this kind of result suggests that the limit joining $\l'$ of Theorem~\ref{thm:main} not only exists, but exhibits some rigidity over different possible initial joinings $\l$, since $\l'$ must exhibit these various instances of relative independence.  Once again there is a superficial analogy here with the study of unipotent flows on homogeneous spaces, where a central theme is the classification of all possible invariant measures and the rigidity that such a classification entails; but once again, I do not know whether this points to any deeper connexions between that setting and ours. \fin

\section{Proof of the main theorem}\label{sec:general-k}

We can now complete the proof of Theorem~\ref{thm:bigmain}. The general case is handled by a `spiral' PET induction on the tuple $(\phi_1,\ldots,\phi_k)$: for each such tuple we will show that
\begin{eqnarray*}
&&(\hbox{assertions (1,2,3) for $(\t{\phi}_2,\ldots,\t{\phi}_k)$})\\
&&\quad\quad\quad \Rightarrow (\hbox{assertion (1) for $(\phi_1,\phi_2,\ldots,\phi_k)$})\\
&&\quad\quad\quad\quad\quad\quad \Rightarrow (\hbox{assertion (2) for $(\phi_1,\phi_2,\ldots,\phi_k)$})\\
&&\quad\quad\quad\quad\quad\quad\quad\quad\quad \Rightarrow (\hbox{assertion (3) for $(\phi_1,\phi_2,\ldots,\phi_k)$}),
\end{eqnarray*}
at which point the induction closes on itself.

We retain the assumptions from the start of Section~\ref{sec:char-factor}. Proposition~\ref{prop:vdC-appn} gives the purchase needed to complete our induction.  Let the class $\sfC$ and group extension $q:\t{G}\to G$ be as in the preceding section.  In analysing the family of averages
\[A^\l_T(f_1,f_2,\ldots,f_k),\]
Proposition~\ref{prop:vdC-appn} allows us to assume that $f_1$ is measurable with respect to the factor $\sfC\S_1$, or equivalently that $\bfX_1$ is itself a system with the property that the $\t{G}$-system $\bfX_1^{q(\cdot)}$ is a factor of a member of the class
\[\sfC_0^{\langle\img\t{\phi}_1\rangle}\vee\bigvee_{j=2}^k\sfC_0^{\langle\img\t{\phi}_j\t{\phi}_1^{-1}\rangle}.\]
From this point a careful re-arrangement gives a reduction to the conclusions of Theorem~\ref{thm:main2} for the group $\t{G}$ and family $(\t{\phi}_2,\ldots,\t{\phi}_k)$, which is isomorphic to $(\phi_2,\ldots,\phi_k)$ and hence precedes $(\phi_1,\phi_2,\ldots,\phi_k)$ in the PET ordering (see Lemma~\ref{lem:PET-calcns}).  Note that this holds in spite of our ascent from $G$ to $\t{G}$, because we have now removed $\t{\phi}_1$ from the picture altogether.

In order to set up the necessary re-arrangement, assume that $\bfX_1 = \sfC\bfX_1$.  By the definition of $\sfC$ there are a system $\t{\bfX} \in \sfC_0^{\langle\img\t{\phi}_1\rangle}\vee\bigvee_{j=2}^k\sfC_0^{\langle\img\t{\phi}_j\t{\phi}_1^{-1}\rangle}$
and a factor map $\pi:\t{\bfX}\to \bfX_1^{q(\cdot)}$.

Now let $\t{\bfX}_1 := \t{\bfX}$ and $\t{\bfX}_i := \bfX_i^{q(\cdot)}$ for any $i \neq 1$, and choose any lift of $\l$ to a joining $\t{\l}$ of the $\t{\bfX}_i$ (for instance, one could use the relatively independent product over $\l$).  For each $i \neq 1$ consider the factor $\t{\S}_1^{\langle \img\t{\phi}_1\t{\phi}_i^{-1}\rangle} \leq \t{\S}_1$, and let
\[\zeta_i:\t{X}_1 \to Z_i\]
be a factor map of standard Borel $\t{G}$-space which generates this factor.  These may be realized as factors of the joining $\t{\l}$ through the coordinate projection $\prod_i\t{X}_i\to \t{X}_1$. Crucially, by enlarging each of the systems $\t{\bfX}_i$ for $i \neq 1$, we can arrange that under $\t{\l}$ each of these factor maps to $Z_i$ is also virtually measurable with respect to the $\t{X}_i$-coordinate, as well as the $\t{X}_1$-coordinate. To this end, for each $i\neq 1$ consider the composition \[\t{X}_0\times\cdots\times \t{X}_k\ \ \stackrel{\rm{coord.}\,\rm{proj.}}{\to}\ \ \t{X}_1\times \t{X}_i\ \stackrel{\zeta_i\times \id}{\to}\ Z_i\times \t{X}_i.\]
Since this composition respects the $\t{G}$-actions, it
defines a joining of $\bfZ_i$ with $\t{\bfX}_i$, which we denote by $\bfY_i = (Y_i,\Phi_i,\nu_i,v_i)$. Let $\eta_i:Y_i\to \t{X}_i$ be the second coordinate projection.

Thus we have constructed a collection of factorizations
\begin{center}
$\phantom{i}$\xymatrix{(\t{X}_0\times \cdots\times \t{X}_k,\t{\S}_0\otimes \cdots \otimes \t{\S}_k,\t{\l},\t{u}_\Delta)\ar[dr]\ar^-{\rm{coord.}\,\rm{proj.}}[rr]&& (\t{X}_i,\t{\S}_i,\t{\mu}_i,\t{u}_i)\\
& (Y_i,\Phi_i,\nu_i,v_i)\ar[ur]
}
\end{center}
for each $i\in \{0,2,3,\ldots,k\}$. Putting these together with the coordinate projection $\t{X}_0\times \cdots\times \t{X}_k \to \t{X}_1$ therefore gives a measure-theoretic isomorphism
\begin{multline*}
(\t{X}_0\times \cdots\times \t{X}_k,\t{\S}_0\otimes \cdots \otimes \t{\S}_k,\t{\l},\t{u}_\Delta)\\
\stackrel{\cong}{\to}(Y_0\times \t{X}_1\times Y_2\times \cdots \times Y_k,\Phi_0\otimes \t{\S}_1\otimes \Phi_2\otimes \cdots \otimes \Phi_k,\theta,v_\Delta)
\end{multline*}
for some joining $\theta$ of $\t{G}$-systems.

In addition, this construction guarantees that the factor maps
\[\t{X}_0\times \cdots\times \t{X}_k\ \ \stackrel{\rm{coord.}\ \rm{proj.}}{\to}\ \ \t{X}_1\stackrel{\zeta_i}{\to} Z_i\]
and
\[\t{X}_0\times \cdots\times \t{X}_k \to  Y_i\ \  \stackrel{\rm{coord.}\,\rm{proj.}}{\to}\ \  Z_i\]
agree up to $\l'$-negligible sets.  Therefore any $h \in L^\infty(\t{\mu}_1)$ which is measurable with respect to $\t{\S}_1^{\langle \img\t{\phi}_1\t{\phi}_i^{-1}\rangle}$ (equivalently, which is invariant under $\langle \img\t{\phi}_1\t{\phi}_i^{-1}\rangle$, with the convention that $\t{\phi}_0 \equiv e$) has an essentially unique counterpart $h' \in L^\infty(\nu_i)$ which lifts to the same function on $\t{X}_0\times \cdots\times \t{X}_k$ up to $\t{\l}$-negligible sets, and which is invariant under the same subgroup of $\t{G}$.

\begin{lem}\label{lem:re-arrange}
In the situation described above, consider the averaging operators associated to the lifted family of polynomial maps $\t{\phi}_i:\bbR\times \bbR^r\to \t{G}$.  Suppose that $f_1 \in L^\infty(\t{\mu}_1)$ is a function of the special form
\[g\cdot h_2\cdot \cdots \cdot h_k,\]
where $g \in L^\infty(\t{\mu}_1)$ is invariant under $\langle \img \t{\phi}_1 \rangle$ and each $h_i$ is invariant under $\langle \img \t{\phi}_1 \t{\phi}_i^{-1} \rangle$. Then for any other functions $f_i \in L^\infty(\t{\mu}_i)$ for $i \neq 1$ one has
\[A^{\t{\l}}_T(f_1,f_2,\ldots,f_k) = \sfE\big(g'\cdot A^\theta_T(1,\,h'_2(f_2\circ\eta_2),\,\ldots,\,h'_k(f_k\circ\eta_k))\,\big|\,\eta_0\big)\]
(recalling that $A^{\t{\l}}_T$ has range in $L^\infty(\t{\mu}_0)$, while $A^\theta_T$ has range in $L^\infty(\t{\nu}_0)$), where $g'$ and $h'_i$ are the counterparts of $g$ and $h_i$ introduced above.
\end{lem}

\textbf{Proof}\quad By the definition of $A^{\t{\l}}_T$ and $A^\theta_T$ this follows from the analogous calculation at the level of joinings.  For the joinings $\t{\l}$ and $\theta$, the above isomorphism gives
\begin{eqnarray*}
&&\int_0^T\int_{\prod_i\t{X}_i} f_0\otimes (f_1\circ \t{u}_1^{\t{\phi}_1(t,h)})\otimes \cdots\otimes (f_k\circ \t{u}_k^{\t{\phi}_k(t,h)})\,\d\t{\l}\,\d t\\
&&= \barint_0^T\int_{Y_0\times \t{X}_1\times Y_2\times \cdots\times Y_k} (f_0\circ \eta_0)\otimes (f_1\circ \t{u}_1^{\t{\phi}_1(t,h)})\\
&&\quad \quad \quad \quad \quad \quad \quad \quad \quad \quad 
\otimes (f_2\circ \eta_2\circ v_2^{\t{\phi}_2(t,h)})\otimes\cdots\otimes (f_k\circ\eta_k\circ v_k^{\t{\phi}_k(t,h)})\,\d\theta\,\d t.
\end{eqnarray*}

On the other hand, our assumptions on the structure of $f_1$ imply that
\[g\circ \t{u}_1^{\t{\phi}_1(t,h)} = g\quad \hbox{and}\quad h_i\circ\t{u}_1^{\t{\phi}_1(t,h)} = h_i\circ\t{u}_1^{\t{\phi}_i(t,h)},\ i=2,3,\ldots,k,\]
for all $(t,h)$.  Also, the counterparts $g' \in L^\infty(\nu_0)$ and $h'_i \in L^\infty(\nu_i)$ for $i\geq 2$ satisfy
\[g'(y_0) = g(\t{x}_1)\quad \hbox{and}\quad h'_i(y_i) = h_i(\t{x}_1)\]
for $\theta$-almost every $(y_0,\t{x}_1,y_2,\ldots,y_k)$.  The above integral with respect to $\theta$ may therefore be re-written as
\begin{multline*}
\barint_0^T\int_{Y_0\times \t{X}_1\times Y_2\times \cdots\times Y_k} (g'(f_0\circ \eta_0))\otimes 1_{\t{X}_1} \otimes  \big((h'_2(f_2\circ \eta_2))\circ v_2^{\t{\phi}_2(t,h)}\big) \otimes\\
\cdots\otimes \big((h'_k(f_k\circ\eta_k))\circ v_k^{\t{\phi}_k(t,h)}\big)\,\d\theta\,\d t.
\end{multline*}
Regarded as a linear functional applied to $f_0$, this is integration against
\[\sfE\big(g\cdot A^\theta_T(1,\,h_2(f_2\circ\eta_2),\,\ldots,\,h_k(f_k\circ\eta_k))\,\big|\,\eta_0\big),\]
as required. \qed

Of course, the importance of the above lemma is that on the right-hand side there is no non-trivial function in the first entry under $A^\theta_T$.  This now leads quite smoothly to a completion of our spiral induction.

\textbf{Proof of Theorem~\ref{thm:bigmain}}\quad In case $k=1$, $M^\l$ extends to a bounded operator $L^2(\mu_1)\to L^2(\mu_0)$ and the desired assertions of convergence and genericity become simply that (i) the average
\[M^\l\Big(\barint_0^T u(\phi_1(t,h)^{-1})^\ast f_1\,\d t\Big)\]
converges to $M^\l P_hf$ with $P_h$ the conditional expectation onto $\S_1^{\langle \img \phi(\cdot,h)\rangle}$, and (ii) this is generically equal to $M^\l P f$ with $P$ the conditional expectation onto $\S_1^{\langle \img \phi\rangle}$.  Both of these assertions follow at once from Proposition~\ref{prop:k=1}.

It remains to handle the inductive step in case $k\geq 2$. Assume that properties (1--3) have already been proved for all tuples preceding $\F$ in the PET ordering.  We will deduce those properties for $\F$ in order.

\quad\textbf{Property (1)}\quad In this step, by fixing one $h$ throughout the proof and replacing $G$ with its subgroup $\langle\img\phi_1(\cdot,h) \cup \cdots \cup \img \phi_k(\cdot,h)\rangle$ if necessary, we may assume that each $\phi_i$ is a function of $t$ alone, and hence that $r=0$. With this agreed, let $q:\t{G}\to G$ and the class $\sfC$ be constructed as before using this new group and tuple of maps.

By re-ordering $\F$ if necessary we may also assume $\phi_1$ is a pivot.  In this case, by Proposition~\ref{prop:vdC-appn} it suffices to show that the averages $A^\l_T(f_1,\ldots,f_k)$ converge when $f_1$ is $\sfC\S_1$-measurable.

Construct the $\t{G}$-systems $\t{\bfX}_i$ and $\bfY_i$ as above. Lifting $f_1$ to $f_1\circ\pi\in L^\infty(\t{\mu}_1)$, on this larger system we know that it can be approximated in $L^2(\t{\mu}_1)$ by finite sums of the form
\[\sum_p g_p\cdot h_{2,p}\cdot \cdots h_{k,p},\]
where $g_p \in L^\infty(\t{\mu}_1)$ is invariant under $\langle \img \t{\phi}_1 \rangle$ and each $h_{i,p} \in L^\infty(\t{\mu}_i)$ is invariant under $\langle \img \t{\phi}_1 \t{\phi}_i^{-1} \rangle$.

Appealing first to the uniform continuity of the operators $A^{\t{\l}}_T$ in each entry separately, and then to the linearity of these operators in the first entry, it therefore suffices to prove convergence of the averages
\[A^{\t{\l}}_T(f_1,\ldots,f_k)\]
whenever $f_1$ is one such product function.  However, this case lands within the hypothesis of the preceding lemma, which converts these into averages of the form
\[\sfE\big(g\cdot A^\theta_T(1,\,h_2(f_2\circ\eta_2),\,\ldots,\,h_k(f_k\circ\eta_k))\,\big|\,\eta_0\big).\]
The norm convergence of these now follows from the norm convergence of the averages $A^\theta_T(1,\,h_2(f_2\circ\eta_2),\,\ldots,\,h_k(f_k\circ\eta_k))$, which is promised by the inductive hypothesis applied to the simpler polynomial family $(\phi_2,\ldots,\phi_k)$.

\quad\textbf{Property (2)}\quad Of course, property (1) already implies convergence of the averaged couplings
\[\barint_0^T (\id_{X_0}\times u_1^{\phi_i(t,h)}\times u_2^{\phi_2(t,h)} \times \cdots \times u_k^{\phi_k(t,h)})_\ast\l\,\d t\]
as $T\to\infty$ to some limit $\l^h$.  We must next show that for any tuple of functions $f_i \in L^\infty(\mu_i)$, the $\l^h$-integrals are the same whether we integrate $f_0\otimes f_1 \otimes \cdots \otimes f_k$ or $(f_0\circ u^{g_0})\otimes (f_1\circ u^{g_1}) \otimes \cdots \otimes (f_k\circ u^{g_k})$ for any
\[(g_0,g_1,\ldots,g_k) \in G^{\Delta (k+1)}\quad\hbox{or}\quad (g_0,g_1,\ldots,g_k) \in\langle \img \vec{\phi}(\cdot,h) \rangle.\]
This will give the invariance of $\l^h$ under the $u_\times$ action of $\langle G^{\Delta (k+1)}\cup\img \vec{\phi}(\cdot,h) \rangle$.

As in the case of property (1), in this step we can fix a choice of $h$ and replace $G$ with the subgroup $G^h := \langle \img\phi_1(\cdot,h)\cup\cdots \cup \img \phi_k(\cdot,h)\rangle$ if necessary, so that we may assume $r = 0$.

Since
\[\sfE(f_1\,|\,\sfC\S_1)\circ u_1^g = \sfE(f_1\circ u_1^g\,|\,\sfC\S_1)\]
for any $g$, by Proposition~\ref{prop:vdC-appn} it again suffices to treat the case when $f_1$ is $(\sfC\S_1)$-measurable.  Now we may consider again the previous construction of the $\t{G}$-systems $\t{\bfX}_i$ and $\bfY_i$ and their joinigs $\l'$ and $\theta$. In these terms we wish to prove that
\begin{multline*}
\int_{\prod_i\t{X}_i}f_0\otimes f_1 \otimes \cdots \otimes f_k\,\d\t{\l}'\\ = \int_{\prod_i\t{X}_i}(f_0\circ \t{u}_0^{\t{g}_0})\otimes (f_1\circ \t{u}_1^{\t{g}_1}) \otimes \cdots \otimes (f_k\circ \t{u}_k^{\t{g}_k})\,\d\t{\l}'
\end{multline*}
for any tuple $f_i \in L^\infty(\t{\mu}_i)$ and any
\[(\t{g}_0,\t{g}_1,\ldots,\t{g}_k) \in \t{G}^{\Delta (k+1)}\quad\hbox{or}\quad (\t{g}_0,\t{g}_1,\ldots,\t{g}_k) \in\langle \img \vec{\t{\phi}}(\cdot) \rangle,\]
where $\t{\l}'$ is the limit joining obtained by averaging $\t{\l}$.

Arguing again as for property (1), by continuity and multilinearity we may now assume that $f_1$ is of the special form $g\cdot h_2\cdot\cdots \cdot h_k$ assumed by Lemma~\ref{lem:re-arrange}, and so by that lemma it now suffices to prove that
\begin{multline*}
\int_{Y_0\times \t{X}_0\times Y_2\times \cdots\times Y_k} (g'(f_0\circ \eta_0))\otimes 1\otimes \cdots \otimes (h'_k(f_k\circ \eta_k))\,\d\theta'\\ = \int_{Y_0\times \t{X}_0\times Y_2\times \cdots\times Y_k} ((g'(f_0\circ \eta_0))\circ v_0^{\t{g}_0})\otimes 1\otimes \cdots \otimes ((h'_k(f_k\circ \eta_k))\circ v_k^{\t{g}_k})\,\d\theta',
\end{multline*}
where $\theta'$ is the limit joining obtained by averaging $\theta$.  With this re-arrangement the coordinate in $\t{X}_1$ vanishes from the picture, and what remains is just an instance of property (2) for the simpler tuple of polynomial maps $(\t{\phi}_2,\ldots,\t{\phi}_k)$, which is known by induction.

\quad\textbf{Property (3)}\quad Lastly, we must show that there is a Zariski residual set $E \subseteq \bbR^r$ such that for any tuple of functions $f_i$ the limit
\[\lim_{T\to\infty}\int_{X_0}f_0\cdot A^\l_T(f_1,\ldots,f_k)\,\d\mu_0\]
is the same for all $h \in  E$, which will imply that the map $h\mapsto \l^h$ is Zariski generically constant (and hence, by property (2), that this generic value must be invariant under the whole of $\langle G^{\Delta (k+1)}\cup\img \vec{\phi} \rangle$). In this step, of course, we may not restrict to a single value of $h$. 

Clearly it suffices to prove this $h$-independence for functions $f_i$ drawn from countable ${\|\cdot\|}_2$-dense subsets of $L^\infty(\mu_i)$, and since a countable intersection of Zariski generic sets is Zariski generic we may therefore look for such a 
Zariski generic set for just a single tuple of functions $f_i$.

The full strength of Proposition~\ref{prop:vdC-appn} and our construction above now give a Zariski residual subset $E \subseteq \bbR^r$, extensions of $\t{G}$-systems $\pi:\t{\bfX}_i\to \bfX_i^{q(\cdot)}$ and a joining $\t{\l}$ of $\t{G}$-systems such that
\begin{multline*}
\int_{X_0\times \cdots\times X_k}f_0\otimes \cdots\otimes f_k\,\d \l^h\\ = \lim_{T\to\infty}\int_{\t{X}_0}(f_0\circ \pi_0)\cdot A^{\t{\l}}_T(\sfE(f_1\circ\pi\,|\,\L),f_2\circ\pi_2,\ldots,f_k\circ \pi_k)\,\d\t{\mu}_0
\end{multline*}
for all $h \in E$, where now
\[\L := \t{\S}_1^{\langle \img \t{\phi}_1\rangle}\vee \bigvee_{i=2}^k \t{\S}_1^{\langle \img \t{\phi}_1\t{\phi}_i^{-1}\rangle}.\]
Clearly it suffices to show that the desired $h$-independence holds on some further Zariski residual subset of $E$, and now the same manipulations as above give a reduction of this to a proof that the limits
\[\lim_{T\to \infty}\int_{Y_0}(g_0'(f_0\circ \eta_0))\cdot A^\theta_T\big(1,h_2'(f_2\circ\eta_2),\ldots,h_k'(f_k\circ\eta_k)\big)\,\d\nu_0\]
are independent of $h$ on some Zariski residual set, where $\theta$ and the $Y_i$ have been constructed from $\l$ and the $\t{X}_i$ as previously.  The dependence on $h$ in this expression is all in the off-diagonal polynomial trajectory that appears in the average $A^\theta_T$.  Once again, the fact that this limit is generically constant now follows from the inductive hypothesis applied to the family $(\t{\phi}_2,\ldots,\t{\phi}_k)$, and so the proof is complete. \qed

\section{Further questions}\label{sec:further-ques}

\subsection{Other questions in continuous time}

Theorems~\ref{thm:main} and~\ref{thm:main2} suggest many possible extensions involving different kinds of averaging, just as for any other equidistribution phenomenon.  The following paragraphs contain a sample of these possibilities.

First, given another connected nilpotent group $G'$, one could ask more generally about polynomial maps $\phi_i:G'\to G$ and the resulting off-diagonal averages along a F\o lner sequence of subsets $F_N \subseteq G'$.  Do these always converge as in our main theorems?  This seems likely, and I suspect that the methods of proof above can provide significant insight into this question, but it may be tricky to set up the right generalization of PET induction.

A little more abstractly, the off-diagonal polynomial trajectory
\[\{(\phi_1(t),\phi_2(t),\ldots,\phi_k(t):\ t\in \bbR\}\]
is a semi-algebraic subset of $G^k$ in the sense of real algebraic geometry (see, for instance, Bochnak, Coste and Roy~\cite{BocCosRoy98}).  Could it be that convergence as in Theorems~\ref{thm:main} or~\ref{thm:main2} holds along the intersections of increasingly large balls with any semi-algebraic subset $V \subset G^k$, endowed with a suitable surface-area measure?

A more challenging question concerns the assumption that $G$ be nilpotent.  Do Theorems~\ref{thm:main} or~\ref{thm:main2} still hold if we assume only that $G$ is an arbitrary connected and simply connected Lie group?  This is probably too much to ask, but some progress may be possible, for instance, if each $\phi_i$ has image lying within a unipotent subgroup of $G$.  This seems a natural setting to investigate in view of Ratner's Theorems giving equidistribution and measure rigidity for unipotent flows on homogeneous spaces~\cite{Rat90-a,Rat90-b,Rat91-a,Rat91-b}, and Shah's extension of these results to averages over regular algebraic maps~\cite{Sha94}.

However, as remarked in the Introduction, the methods used to study homogeneous space flows are very different (and mostly much more delicate) from those explored in this paper.  Shah's analysis of regular algebraic maps proceeds by first obtaining the invariance of a weak limit measure under some unipotent subgroup and then using the resulting structure promised by Ratner's Theorems, whereas it is an essential feature of our inductive proof of Theorem~\ref{thm:bigmain} that the cases of homomorphisms $\phi_i$ and of more general polynomial maps must be treated together.

To illustrate more concretely some of the difficulties posed by non-nilpotent groups, consider the functional averages
\[\barint_0^T (f_1\circ u_1^t)(f_2\circ u_2^t)\,\d t\]
for a jointly measurable probability-preserving system $(X,\S,\mu,u)$ for $G = \rm{SL}_2(\bbR)$ and with $u_1,u_2:\bbR\to \rm{SL}_2(\bbR)$ parametrizing the upper- and lower-triangular subgroups respectively.  (These averages are easily expressed in terms of the natural analog of Theorem~\ref{thm:main2}.)  If we assume that these averages do not tend to $0$ for some choice of $f_1,f_2 \in L^\infty(\mu)$, then the van der Corput estimate and a re-arrangement give also
\[\barint_0^S\barint_0^T \int_X f_1\cdot (f_1\circ u_1^s)\cdot ((f_2\cdot (f_2\circ u_2^s))\circ u_2^tu_1^{-t})\,\d\mu\,\d t\,\d s \,\,\not\!\!\to 0\]
as $T\to\infty$ and then $S\to \infty$.  In order to use this, we need some information about the averages along the trajectory $t\mapsto u_2^tu_1^{-t}$ in $G$.  This is certainly a polynomial map in the sense of real algebraic geometry, but not in the sense of Definition~\ref{dfn:poly}, so further differencing does not seem to lead to a simplification of the problem.  I have not examined in detail what other arguments (for example, using the representation theory of $\rm{SL}_2(\bbR)$) might be brought to bear here, since this is only a very special case: it simply serves to illustrate that the method of PET induction cannot be applied so na\"\i vely in this setting.

Finally, linked to the study of convergence and equidistribution is the problem of describing the limit joinings $\l'$.  Some information on their possible structure is contained in the proof of Proposition~\ref{prop:vdC-appn} above, as remarked after that proposition, but it would be interesting to know whether they can be classified more precisely, possibly after extending each $\bfX_i$ to a suitably-sated extension.  A discussion of related issues in the setting of $\bbZ^d$-actions can be found in~\cite{Aus--thesis}.

\subsection{Discrete actions}\label{subs:compare-discrete}

Most past interest in the kind of off-diagonal average appearing in Theorem~\ref{thm:main} has focused on actions of discrete groups.  Suppose that $\G$ is a discrete nilpotent group, $\phi_1,\phi_2,\ldots,\phi_k:\bbZ\to\G$ are polynomial maps (according to the obvious relative of Definition~\ref{dfn:poly}), $\bfX_i = (X_i,\S_i,\mu_i,T_i)$ are probability-preserving $\G$-systems for $1 \leq i \leq k$ and $\l$ is a joining of the 
systems $\bfX_i$.  Much recent work has been directed towards understanding whether the off-diagonal averages
\[\frac{1}{N}\sum_{n=1}^N (T_1^{\phi_1(n)}\times \cdots\times T_k^{\phi_k(n)})_\ast\l\]
converge to some limit joining as $N\to\infty$, or whether the associated functional averages converge.  Several partial results have appeared, and at the time of this writing Miguel Walsh has just settled the general case in his preprint~\cite{Wal11}.

Walsh's approach does not use heavy ergodic-theoretic machinery.  It relies on reformulating the problem of norm convergence for the functional averages into a problem asking for some `quantitative' guarantee that one can find long intervals of times $N$ in which those averages are all close in $\|\cdot\|_2$.  This new assertion can then be proved by a clever induction on the tuple of polynomial maps $(\phi_1,\ldots,\phi_k)$, which is apparently different from Bergelson's PET induction.

In making this reformulation, Walsh uses ideas that have some precedent in Tao's proof of convergence when $\G = \bbZ^d$ and all the $\phi_i$ are linear~(\cite{Tao08(nonconv)}).  Some of these ideas lie outside more traditional ergodic-theoretic approaches to this class of questions (such as the present paper), and they have the consequence that very little can be gleaned about the structure of the limits (functions or joinings).  Therefore it would still be of interest to see a proof that gives some additional information, similar to our Theorem~\ref{thm:bigmain} or to the earlier, even more precise results of~\cite{HosKra05} or~\cite{Zie07} in the case of discrete powers of a single transformation.  We finish with an informal discussion of the difficulties that face any attempt to adapt the arguments of the preceding sections to the setting of discrete $\G$.

The first and most obvious difficulty is that if these averaged couplings do converge to some limit $\l'$, it need not be invariant under the off-diagonal subgroup
\[\langle \img (\phi_1,\ldots,\phi_k) \rangle \leq \G^k.\]

Indeed, let $\G = \bbZ$, let $\phi_1 \equiv 0$ and $\phi_2(n) := n^2$, and let $\bfX_1 = \bfX_2$ be the system given by the generator rotation on $\bbZ/4\bbZ$.  Since all square numbers are congruent to either $0$ or $1 \!\!\mod 4$, it is easily computed that the limit obtained by averaging the diagonal joining $\l$ is simply
\[\frac{1}{2}\l + \frac{1}{2}(\id\times T)_\ast\l,\]
which is not $(\id\times T)$-invariant.

Of course, this is a trivial example, but it is not clear whether this kind of arithmetic system, appearing as a factor of more general systems $\bfX_i$, is the only possible obstruction to the desired extra invariance of the limit joining.

While this example bears only on the possible symmetries of the limit joining, in the continuous-time setting those symmetries play a crucial r\^ole in the proof of Proposition~\ref{prop:vdC-appn} above, and so the whole method of proof we have used in this paper may need substantial modification before it can give convergence results in the discrete-time world.

A second difficulty worth remarking is the absence of any useful replacement for the notion of Zariski genericity in the discrete-time setting.  Of course, Corollary~\ref{cor:generically-const-fpspace} is still true for discrete group actions: the problem is that it tells us nothing, because these groups are themselves countable.

It might be worth exploring a more subtle appeal to the reasoning of Corollary~\ref{cor:rel-ind-over-common} in place of Corollary~\ref{cor:generically-const-fpspace}.  The statement of Corollary~\ref{cor:rel-ind-over-common} is also still true for discrete groups provided the subgroups $H_1$ and $H_2$ are both normal in $\langle H_1\cup H_2\rangle$.  One possibility might begin as follows. If $\frH_1$, $\frH_2$, \ldots, is a sequence of closed subspaces of a Hilbert space $\frH$, any two of which are relatively orthogonal over some common further subspace $\frak{K}$, and if in addition $x \in \frH$ is such that $\inf_n\|P_nx\| > 0$ with $P_n$ the orthoprojection onto $\frH_n$, then $x$ also has a nonzero projection onto $\frak{K}$ (for otherwise the $P_nx$ would be an infinite sequence of mutually orthogonal projections of a single vector, all of them large, contradicting Bessel's Inequality).

Structure like this has previously been identified within orthogonal representations of a finitely generated nilpotent group by Leibman~\cite{Lei00}.  Using this reasoning, for example, one can show that if
\[\G = \langle a,b\,|\,[a,b] =: c\ \hbox{is central}\rangle\]
is the discrete Heisenberg group and $T:\G\actson (X,\S,\mu)$ is any action of it, then the $\s$-subalgebras
\[\S^{\langle a \rangle} := \{A \in \S:\ \mu(T^aA\triangle A) = 0\}\]
and $\S^{\langle b\rangle}$ are relatively independent over the fully invariant factor $\S^T$, even though in this discrete setting it can happen that $\S^{\langle a\rangle} \neq \S^{\langle a\rangle^\rm{n}}$ and $\S^{\langle a\rangle}$ is not globally $T$-invariant.  This follows because a judicious appeal to the discrete version of Corollary~\ref{cor:rel-ind-over-common} implies that the $\s$-algebras
\[\S^{b^k\langle a\rangle b^{-k}},\quad k\in\bbZ,\]
are all relatively independent over $\S^{\langle a,c\rangle}$, where $\langle a,c\rangle$ \emph{is} normal in $\G$.  If now $f$ and $g$ are $T^a$- and $T^b$-invariant respectively, then applying $T^b$ gives
\[\int f \cdot \sfE(g\,|\,\S^{\langle a\rangle})\, \d\mu = \int (f \cdot \sfE(g\,|\,\S^{\langle a\rangle}))\circ T^{b^k}\, \d\mu = \int (f\circ T^{b^k})\cdot(\sfE(g\,|\,\S^{b^{-k}\langle a\rangle b^k})\, \d\mu.\]
Therefore the non-vanishing of this integral implies that $g$ actually has uniformly nonzero conditional expectation onto every $\S^{b^k\langle a\rangle b^{-k}}$.  Hence by the argument sketched above, it must actually have nonzero conditional expectation onto $\S^{\langle a,c\rangle}$, and similarly $f$ must have nonzero conditional expectation onto $\S^{\langle b, c\rangle}$.  These two $\s$-algebras are now globally $T$-invariant and relatively independent over $\S^T$, so putting this together shows that $\S^{\langle a\rangle}$ and $\S^{\langle b\rangle}$ are themselves relatively independent over $\S^T$.

In order to use a similar idea to study off-diagonal or multiple averages, one might, for instance, try to prove a discrete analog of Proposition~\ref{prop:vdC-appn} according to which the characteristic factors $\L_h$ obtained depending on $h$ are not mostly equal to each other, but are all relatively orthogonal over some common smaller $\s$-algebra $\L'$. Then it might be possible to replace $\L^h$ with $\L'$ in subsequent arguments and gain more purchase on the asymptotic behaviour of our averages as a result.  However, I do not have a precise statement to formulate based on this speculation.

\textbf{Acknowledgements}\quad This work was supported by a research fellowship from the Clay Mathematics Institute.  Much of it was carried out during a visit to the Isaac Newton Institute for the Mathematical Sciences. \fin

\appendix

\section{A continuous-time van der Corput
estimate}\label{app:cts-time-vdC}

We recall here for completeness a continuous-time variant of the
classical van der Corput estimate for bounded Hilbert-space-valued
sequences.  The discrete-time version can be found in Section 1
of~\cite{FurWei96}, and a continuous-time version in Appendix B of Potts~\cite{Pot09}.

\begin{lem}\label{lem:vdC}
If $u:[0,\infty)\to \frH$ is a bounded strongly measurable map into
a Hilbert space, then vector-valued non-convergence
\[\barint_0^T u(t)\,\d t \,\,\not\!\!\to 0\quad\quad\hbox{as}\ T\to\infty\]
implies the scalar-valued non-convergence
\[\barint_0^S\barint_0^T\langle u(t+s),u(t)\rangle\,\d t\,\d s \,\,\not\!\!\to 0\quad\quad\hbox{as}\ T\to\infty\ \hbox{and then}\ S\to\infty.\]
\qed
\end{lem}

\bibliographystyle{abbrv}
\bibliography{cts_nilp5}

\parskip 0pt

\noindent \small{Mathematics Department, Brown University,}

\noindent \small{Box 1917, 151 Thayer Street,}

\noindent \small{Providence, RI 02912, USA}

\noindent \small{\texttt{timaustin@math.brown.edu},}

\noindent \small{\texttt{www.math.brown.edu/$\sim$timaustin}}

\end{document}